\documentclass[10pt,reqno]{amsart}
\usepackage{mathptmx}
\usepackage{helvet}
\usepackage{courier}
\usepackage{amssymb}
\usepackage{microtype}
\usepackage{hyperref}
\usepackage{mathtools,amsthm}
\usepackage{enumerate}
\usepackage{natbib}
\usepackage{bm}
\usepackage{paralist}
\usepackage{tikz}
\newcommand{\heads}{\mathit{heads}}
\newcommand{\tails}{\mathit{tails}}

\newcommand{\reals}{\mathbb{R}}
\newcommand{\nats}{\mathbb{N}}
\newcommand{\booleans}{\mathbb{B}}
\newcommand{\unit}{[0,1]}
\newcommand{\rv}{X}
\newcommand{\brv}{\mathbf{\rv}}
\newcommand{\arv}{Y}
\newcommand{\values}{\mathcal{X}}
\newcommand{\avalues}{\mathcal{Y}}

\newcommand{\pr}{P}
\newcommand{\lpr}{{\underline{\pr}}}
\newcommand{\upr}{{\overline{\pr}}}
\newcommand{\apr}{Q}
\newcommand{\alpr}{{\underline{\apr}}}

\newcommand{\yapr}{R}
\newcommand{\yalpr}{{\underline{\yapr}}}

\newcommand{\nex}{E}
\newcommand{\lnex}{{\underline{\nex}}}
\newcommand{\unex}{{\overline{\nex}}}
\newcommand{\cylext}[1]{{\widetilde{#1}}}
\newcommand{\muhy}{{\mathit{MuHy}}}
\newcommand{\hy}{{\mathit{Hy}}}
\newcommand{\mult}{{\mathit{Mn}}}
\newcommand{\cmult}{{\mathit{CoMn}}}
\newcommand{\cbin}{{\mathit{CoBi}}}
\newcommand{\gambles}{\mathcal{L}}
\newcommand{\conts}{\mathcal{C}}

\newcommand{\linprevs}{\mathbb{P}}
\newcommand{\domain}{\mathcal{K}}
\newcommand{\adomain}{\mathcal{H}}
\newcommand{\atoms}{\mathcal{A}}

\newcommand{\permuts}{\mathcal{P}}

\newcommand{\btheta}{{\bm{\theta}}}
\newcommand{\bmu}{{\bm{\mu}}}
\newcommand{\simplex}{\Sigma}
\newcommand{\poly}{\mathcal{V}}
\newcommand{\counts}{\mathcal{N}}
\newcommand{\sample}[1]{{\mathbf{#1}}}
\newcommand{\cnt}[1]{{\mathbf{#1}}}
\newcommand{\cntf}{T}
\newcommand{\bcntf}{\mathbf{\cntf}}
\newcommand{\freqf}{F}
\newcommand{\bfreqf}{\mathbf{\freqf}}
\newcommand{\set}[2]{\left\{#1\colon#2\right\}}

\newcommand{\tuple}[2]{({#1}_1,\dots,{#1}_{#2})}
\newcommand{\ntuple}[2]{{#1}_1,\dots,{#1}_{#2}}
\newcommand{\abs}[1]{\lvert#1\rvert}

\newcommand{\solp}{\mathcal{M}}

\DeclareMathOperator{\lproj}{\underline{proj}}

\newcommand{\citegen}[1]{\citeauthor{#1}'s \citeyearpar{#1}}
\newcommand{\citegensec}[2]{\citeauthor{#2}'s \citeyearpar[#1]{#2}}
\theoremstyle{plain}
\newtheorem{theorem}{Theorem}
\newtheorem{proposition}[theorem]{Proposition}

\newtheorem{lemma}[theorem]{Lemma}
\theoremstyle{remark}
\newtheorem*{running}{Running example}

\theoremstyle{definition}

\newcounter{saveenumi}
\begin{document}
\title{Exchangeable lower previsions}
\author{Gert de Cooman \and Erik Quaeghebeur}
\address{Ghent University, SYSTeMS Research Group,
  Technologiepark--Zwijnaarde 914, 9052 Zwijnaarde, Belgium}
\email{gert.decooman@ugent.be, erik.quaeghebuer@ugent.be}
\author{Enrique Miranda}
\address{Rey Juan Carlos University, Dept.~of Statistics and
Operations Research. C-Tulip{\'a}n, s/n, 28933, M{\'o}stoles, Spain}
\email{enrique.miranda@urjc.es}

\begin{abstract}
  We extend \citegen{finetti1937} notion of exchangeability to finite and
  countable sequences of variables, when a subject's beliefs about them are
  modelled using coherent lower previsions rather than (linear) previsions.
  We prove representation theorems in both the finite and the countable case,
  in terms of sampling without and with replacement, respectively. We also
  establish a convergence result for sample means of exchangeable sequences.
  Finally, we study and solve the problem of exchangeable natural extension:
  how to find the most conservative (point-wise smallest) coherent and
  exchangeable lower prevision that dominates a given lower prevision.
\end{abstract}
\keywords{Exchangeability, lower prevision, Representation Theorem,
Bernstein polynomials, convergence in distribution, exchangeable
natural extension, sampling without replacement, multinomial
sampling, imprecise probability, coherence.} \maketitle

\section{Introduction}\label{sec:intro}
This paper deals with belief models for both finite and countable sequences of
exchangeable random variables taking a finite number of values. When such
sequences of random variables are assumed to be exchangeable, this more or
less means that the specific order in which they are observed is deemed
irrelevant.
\par
The first detailed study of exchangeability was made by \citet{finetti1937}
(with the terminology of `equivalent' events). He proved the now famous
Representation Theorem, which is often interpreted as stating that a sequence
of random variables is exchangeable if it is conditionally independent and
identically distributed (IID).\footnote{See \citet[Section~11.4]{finetti1975};
  and \citet{cifarelli1996} for an overview of de Finetti's work.}  Other
important work on exchangeability was done by, amongst many others,
\citet{hewitt1955}, \citet{heath1976}, \citet{diaconis1980} and, in the
context of the behavioural theory of imprecise probabilities that we are going
to consider here, by \citet{walley1991}. We refer to
\citet{kallenberg2002,kallenberg2005} for modern, measure-theoretic
discussions of exchangeability.
\par
One of the reasons why exchangeability is deemed important, especially by
Bayesians, is that, by virtue of de Finetti's Representation Theorem, an
exchangeable model can be seen as a convex mixture of multinomial models.
This has given some ground \citep{finetti1937,finetti1975,dawid1985} to the
claim that aleatory probabilities and IID processes can be eliminated from
statistics, and that we can restrict ourselves to considering exchangeable
sequences instead.\footnote{For a critical discussion of this claim, see
  \citet[Section~9.5.6]{walley1991}.}
\par
De Finetti presented his study of exchangeability in terms of the behavioural
notion of previsions, or fair prices. The central assumption underlying his
approach is that a subject should be able to specify a fair price $\pr(f)$ for
any risky transaction (which we shall call a \emph{gamble}) $f$
\citep[Chapter~3]{finetti1974}. This is tantamount to requiring that he should
always be willing and able to decide, for any real number $r$, between selling
the gamble $f$ for $r$, or buying it for that price. This may not always be
realistic, and for this reason, it has been suggested that we should
explicitly allow for a subject's indecision, by distinguishing between his
\emph{lower prevision} $\lpr(f)$, which is the supremum price for which he is
willing to buy the gamble $f$, and his \emph{upper prevision} $\upr(f)$, which
is the infimum price for which he is willing to sell $f$. For any real number
$r$ strictly between $\lpr(f)$ and $\upr(f)$, the subject is then not
specifying a choice between selling or buying the gamble $f$ for $r$. Such
lower and upper previsions are also subject to certain rationality or
coherence criteria, in very much the same way as (precise) previsions are on
de~Finetti's account. The resulting \emph{theory of coherent lower
  previsions}, sometimes also called the behavioural theory of imprecise
probabilities, and brilliantly defended by \cite{walley1991}, generalises de
Finetti's behavioural treatment of subjective, epistemic probability, and
tries to make it more realistic by allowing for a subject's indecision. We
give a brief overview of this theory in Section~\ref{sec:distributions}.
\par
Also in this theory, it is interesting to consider what are the consequences
of a subject's \emph{exchangeability assessment}, i.e., that the order in
which we consider a number of random variables is of no consequence. This is
our motivation for studying exchangeable \emph{lower} previsions in this
paper. An assessment of exchangeability will have a clear impact on the
structure of so-called \emph{exchangeable} coherent lower previsions. We shall
show they can be written as a combination of
\begin{inparaenum}[(i)]
\item a coherent (linear) prevision expressing that permutations of
  realisations of such sequences are considered equally likely, and
\item a coherent lower prevision for the `frequency' of occurrence of the
  different values the random variables can take.
\end{inparaenum}
Of course, this is the essence of representation in de Finetti's sense: we
generalise his results to coherent lower previsions.
\par
A subject's probability assessments may be \emph{local}, in the sense that
they concern the probabilities or previsions of specific events or random
variables. Assessments may on the other hand also be \emph{structural}
\citep[see][Chapter~9]{walley1991}, in which case they specify relationships
that should hold between the probabilities or previsions of a number of events
or random variables. One may wonder if (and how) it is possible to combine
local with structural assessments, such as exchangeability. We show that this
is indeed the case, and give a surprisingly simple procedure, called
\emph{exchangeable natural extension}, for finding the point-wise smallest
(most conservative) coherent and exchangeable lower prevision that dominates
the local assessments. As an example, we use our conclusions to take a fresh
look at the old question whether a given exchangeable model for $n$ variables
can be extended to an exchangeable model for $n+k$ variables.
\par
Before we go on, we want to draw attention to a number of distinctive features
of our approach. First of all, the usual proofs of the Representation Theorem,
such as the ones given by \citet{finetti1937}, \citet{heath1976}, or
\citet{kallenberg2005}, do not lend themselves very easily to a generalisation
in terms of coherent lower previsions. In principle it would be possible, at
least in some cases, to start with the versions already known for (precise)
previsions, and to derive their counterparts for lower previsions using
so-called lower envelope theorems (see Section~\ref{sec:distributions} for
more details). This is the method that \citet[Sections~9.5.3
and~9.5.4]{walley1991} suggests. But we have decided to follow a different
route: we derive our results directly for lower previsions, using an approach
based on Bernstein polynomials, and we obtain the ones for previsions as
special cases.  We believe this method to be more elegant and self-contained,
and it certainly has the additional benefit of drawing attention to what we
feel is the essence of de Finetti's Representation Theorem: specifying a
coherent belief model for a countable exchangeable sequence is tantamount to
specifying a coherent (lower) prevision on the linear space of polynomials on
some simplex, and nothing more.
\par
Secondly, we shall focus on, and use the language of, (lower and upper)
previsions for gambles, rather than (lower and upper) probabilities for
events. Our emphasis on prevision or expectation, rather than probability, is
in keeping with \citegen{finetti1974} and \citegen{whittle2000} approach to
probabilistic modelling. But it is not merely a matter of aesthetic
preference: as we shall see, in the behavioural theory of imprecise
probabilities, the language of gambles is much more expressive than that of
events, and we need its full expressive power to derive our results.
\par
The plan of the paper is as follows. In Section~\ref{sec:distributions}, we
introduce a number of results from the theory of coherent lower previsions
necessary to understand the rest of the paper. In
Section~\ref{sec:exchangeable-finite}, we define exchangeability for finite
sequences of random variables, and establish a representation of coherent
exchangeable lower previsions in terms of sampling without replacement. In
Section~\ref{sec:exchangeable-sequences}, we extend the notion of
exchangeability to countable sequences of random variables, and in
Section~\ref{sec:representation-sequences} we generalise de Finetti's
Representation Theorem (in terms of multinomial sampling) to exchangeable
coherent lower previsions.  The results we obtain allow us to develop a limit
law for sample means in Section~\ref{sec:sample-means}.
Section~\ref{sec:exchangeable-natex} deals with exchangeable natural
extension: combining local assessments with exchangeability.  In an appendix,
we have gathered a few useful results about multivariate Bernstein
polynomials.

\section{Lower previsions, random variables and their distributions}
\label{sec:distributions}
In this section, we want to provide a brief summary of ideas, and known as
well as new results from the theory of coherent lower previsions
\citep{walley1991}. This should lead to a better understanding of the
developments in the sections that follow. For results that are mentioned
without proof, proofs can be found in \citet{walley1991}.

\subsection{Epistemic uncertainty models}
Consider a \emph{random variable} $\rv$ that may assume values $x$ in some
non-empty set $\values$. By `random', we mean that a subject is uncertain
about the actual value of the variable $\rv$, i.e., does not know what this
actual value is. But we do assume that the actual value of $\rv$ can be
determined, at least in principle. Thus we may for instance consider tossing a
coin, where $\rv$ is the outcome of the coin toss, and
$\values=\{\heads,\tails\}$. It does not really matter here to distinguish
between a subject's belief before tossing the coin, or after the toss where,
say, the outcome has been kept hidden from the subject. All that matters for
us here is that our subject is in a state of (partial) ignorance because of a
lack of knowledge. The uncertainty models that we are going to describe here
are therefore \emph{epistemic}, rather than physical, probability models.
\par
Our subject may be uncertain about the value of $\rv$, but he may entertain
certain beliefs about it. These beliefs may lead him to engage in certain
risky transactions whose outcome depends on the actual value of $\rv$. We are
going to try and model his beliefs mathematically by zooming in on such risky
transactions. They are captured by the mathematical concept of a \emph{gamble}
on $\values$, which is a bounded map $f$ from $\values$ to the set $\reals$ of
real numbers. A gamble $f$ represents a random reward: if the subject
\emph{accepts} $f$, this means that he is willing to engage in the following
transaction: we determine the actual value~$x$ that~$\rv$ assumes in~$\values$, and then the subject receives the (possibly negative) reward
$f(x)$, expressed in units of some predetermined linear utility. Let us denote
by $\gambles(\values)$ the set of all gambles on $\values$.
\par
\Citet{finetti1974} has proposed to model a subject's beliefs by eliciting his
fair price, or \emph{prevision}, $\pr(f)$ for certain gambles~$f$. This
$\pr(f)$ can be defined as the unique real number~$p$ such that the subject is
willing to buy the gamble $f$ for all prices $s$ (i.e., accept the gamble~$f-s$) and sell $f$ for all prices $t$ (i.e., accept the gamble $t-g$) for all
$s<p<t$. The problem with this approach is that it presupposes that there is
such a real number, or, in other words, that the subject, whatever his beliefs
about $\rv$ are, is willing, for (almost) every real $r$, to make a choice
between buying $f$ for the price $r$, or selling it for that price.

\subsection{Coherent lower previsions and natural
  extension}\label{sec:coherence}
A way to address this problem is to consider a model which allows our subject
to be undecided for some prices~$r$. This is done in \citegen{walley1991}
theory of lower and upper previsions. The \emph{lower prevision} of the gamble
$f$, $\lpr(f)$, is our subject's supremum acceptable buying price for $f$;
similarly, our subject's \emph{upper prevision}, $\upr(f)$, is his infimum
acceptable selling price for~$f$. Hence, he is willing to buy the gamble~$f$
for all prices $t<\lpr(f)$ and sell $f$ for all prices $s>\upr(f)$, but he may
be undecided for prices $\lpr(f)\leq p\leq\upr(f)$.
\par
Since buying the gamble $f$ for a price $t$ is the same as selling the gamble
$-f$ for the price $-t$ [in both cases we accept the gamble $f-t$], the lower
and upper previsions are \emph{conjugate} functions: $\lpr(f)=-\upr(-f)$ for
any gamble $f$. This allows us to concentrate on one of these functions, since
we can immediately derive results for the other. In this paper, we focus
mainly on lower previsions.
\par
If a subject has made assessments about the supremum buying price (lower
prevision) for all gambles in some domain $\domain$, we have to check that
these assessments are consistent with each other.  First of all, we say that
the lower prevision $\lpr$ \emph{avoids sure loss} when
\begin{equation}\label{eq:asl}
  \sup_{x\in\values}
  \left[
    \sum_{k=1}^{n}\lambda_k[f_k(x)-\lpr(f_k)]\right]
  \geq0
\end{equation}
for any natural number $n$, any gambles $f_1$, \dots, $f_n$ in $\domain$ and
any non-negative real numbers $\lambda_1$, \dots, $\lambda_n$. When the
inequality~\eqref{eq:asl} is not satisfied, there is some non-negative
combination of acceptable transactions that results in a transaction that
makes our subject lose utiles, no matter the outcome, and we then say that his
lower prevision $\lpr$ \emph{incurs sure loss}.
\par
More generally, we say that the lower prevision $\lpr$ is \emph{coherent} when
\begin{equation}\label{eq:coherence}
  \sup_{x\in\values}
  \left[
    \sum_{k=1}^{n}\lambda_k[f_k(x)-\lpr(f_k)]-\lambda_0[f_0(x)-\lpr(f_0)]
  \right]
  \geq0
\end{equation}
for any natural number $n$, any gambles $f_0$, \dots, $f_n$ in $\domain$ and
any non-negative real numbers $\lambda_0$, \dots, $\lambda_n$. Coherence means
that our subject's supremum acceptable buying price for a gamble $f$ in the
domain cannot be raised by considering the acceptable transactions implicit in
other gambles.  In particular, it means that $\lpr$ avoids sure loss. We call
an upper prevision coherent if its conjugate lower prevision is.
\par
If a lower prevision $\lpr$ is defined on a linear space of gambles $\domain$,
then the coherence requirement~\eqref{eq:coherence} is equivalent to the
following conditions: for any gambles $f$ and $g$ in $\domain$ and any
non-negative real number $\lambda$, it should hold that:
\begin{enumerate}[(P1)]
\item $\lpr(f)\geq \inf f$ [accepting sure gains];\label{item:asg}
\item $\lpr(\lambda f)=\lambda \lpr(f)$ [non-negative
  homogeneity];\label{item:nonneghom}
\item $\lpr(f+g)\geq \lpr(f)+\lpr(g)$ [super-additivity].\label{item:superadd}
\end{enumerate}
Moreover, a lower prevision on a general domain is coherent if and only if it
can be extended to a coherent lower prevision on some linear space.
\par
A coherent lower prevision that is defined on indicators of events only is
called a coherent \emph{lower probability}. The indicator $I_A$ of an event
$A$ is the $\{0,1\}$-valued gamble given by $I_A(x):=1$ if $x \in A$ and
$I_A(x):=0$ otherwise.
\par
On the other hand, a lower prevision $\lpr$ on some set of gambles $\domain$
that avoids sure loss can always be `corrected' and extended to a coherent
lower prevision on $\gambles(\values)$, in a least-committal manner: the
(point-wise) smallest, and therefore most conservative, coherent lower
prevision on $\gambles(\values)$ that (point-wise) dominates $\lpr$ on
$\domain$, is called the \emph{natural extension} of $\lpr$, and it is given
for all $f$ in $\gambles(\values)$ by
\begin{equation}\label{eq:natex}
  \lnex(f)
  :=\sup\set{
    \inf_{x\in\values}\bigg[f(x)-\sum_{k=1}^n\lambda_k[f_k(x)-\lpr(f_k)]\bigg]}
  {n\geq0,\lambda_k\geq0,f_k\in\domain}.
\end{equation}
The natural extension of $\lpr$ provides the supremum acceptable buying prices
that we can derive for any gamble $f$ taking into account only the buying
prices for the gambles in $\domain$ and the notion of coherence.
Interestingly, $\lpr$ is coherent if and only if it coincides with its natural
extension $\lnex$ on its domain $\domain$, and in that case $\lnex$ is the
point-wise smallest coherent lower prevision that extends $\lpr$ to
$\gambles(\values)$.

\subsection{Linear previsions}
If the lower prevision $\lpr(f)$ and the upper prevision $\upr(f)$ for a
gamble $f$ happen to coincide, then the common value $\pr(f)=\lpr(f)=\upr(f)$
is called the subject's (precise) \emph{prevision} for $f$. Previsions are
fair prices in \citegen{finetti1974} sense. We shall call them \emph{precise}
probability models, and lower previsions will be called \emph{imprecise}.
Specifying a prevision $\pr$ on a domain $\domain$ is tantamount to specifying
both a lower prevision $\lpr$ and an upper prevision $\upr$ on $\domain$ such
that $\lpr(f)=\upr(f)=\pr(f)$. Since then, by conjugacy,
$\upr(f)=-\lpr(-f)=-\pr(-f)$, it is also equivalent to specifying a lower
prevision $\lpr$ on the larger and negation invariant domain
$\domain':=\domain\cup-\domain$, by letting $\lpr(f):=\pr(f)$ if $f\in\domain$
and $\lpr(f):=-\pr(-f)$ if $f\in-\domain$.  This prevision $\pr$ is then
called coherent, or \emph{linear}, if and only if the associated lower
prevision $\lpr$ is coherent, and this is equivalent to the following
condition
\begin{equation*}
  \sup_{x\in\values}
  \left[
    \sum_{k=1}^{n}\lambda_k[f_k(x)-\pr(f_k)]
    -\sum_{\ell=1}^m\mu_\ell[g_\ell(x)-\pr(g_\ell)]
  \right]
  \geq0
\end{equation*}
for any natural numbers $n$ and $m$, any gambles $f_1$, \dots, $f_n$ and
$g_1$, \dots, $g_m$ in $\domain$ and any non-negative real numbers
$\lambda_1$, \dots, $\lambda_n$ and $\mu_1$, \dots, $\mu_m$.
\par
A prevision on the set $\gambles(\values)$ of all gambles is linear if and
only if it is a positive ($f\geq0\Rightarrow\pr(f)\geq0$) and normed
($\pr(1)=1$) real linear functional. A prevision on a general domain is linear
if and only if it can be extended to a linear prevision on all gambles. We
shall denote by $\linprevs(\values)$ the set of all linear previsions on
$\gambles(\values)$.
\par
The restriction of a linear prevision $\pr$ on $\gambles(\values)$ to the set
$\wp(\values)$ of (indicators of) all events, is a finitely additive
probability. Conversely, a finitely additive probability on $\wp(\values)$ has
a unique extension (namely, its natural extension as a coherent lower
probability) to a linear prevision on $\gambles(\values)$. In this sense, such
linear previsions and finitely additive probabilities can be considered
equivalent: for precise probability models, the language of events is as
expressive as that of gambles.
\par
A linear prevision that is defined on indicators of events only, and therefore
called a coherent probability, is always the restriction of some finitely
additive probability.
\par
There is an interesting link between precise and imprecise probability models,
expressed through the following so-called \emph{lower envelope theorem}: A
lower prevision $\lpr$ on some domain $\domain$ is coherent if and only if it
is the \emph{lower envelope} of some set of linear previsions, and in
particular of the convex set $\solp(\lpr)$ of all linear previsions that dominate it:
for all $f$ in $\domain$,
\begin{equation*}
  \lpr(f)=\inf\set{\pr(f)}{\pr\in\solp(\lpr)},
\end{equation*}
where $\solp(\lpr):=\set{\pr\in\linprevs(\values)}{(\forall
  f\in\domain)(\lpr(f)\geq\pr(f))}$.  We can also use the set $\solp(\lpr)$ to
calculate the natural extension of $\lpr$: for any gamble $f$ on $\values$, we
have that
\begin{equation*}
  \lnex(f):=\inf\set{\pr(f)}{\pr\in\solp(\lpr)}.
\end{equation*}
\par
If we have a coherent lower probability defined on some set of events, then
there will generally be many (i.e., an infinity of) coherent lower previsions
that extend it to all gambles. In this sense, the language of gambles is
actually \emph{more expressive} than that of events when we are considering
lower rather than precise previsions. As already signalled in the
Introduction, this is the main reason why in the following sections, we shall
formulate our study of exchangeable lower previsions in terms of gambles and
lower previsions rather than events and lower probabilities.

\subsection{Important consequences of coherence}
Let us list a few consequences of coherence that we shall have occasion to use
further on. Besides the properties (P\ref{item:asg})--(P\ref{item:superadd})
we have already mentioned that hold when the domain of $\lpr$ is a linear
space, the following properties hold for a coherent lower prevision whenever
the gambles involved belong to its domain:
\begin{enumerate}[(i)]
\item $\lpr$ is \emph{monotone}: if $f\leq g$, then $\lpr(f)\leq\lpr(g)$.
\item $\inf f\leq\lpr(f)\leq\upr(f)\leq\sup f$.
\end{enumerate}
Moreover, coherent lower and upper previsions are continuous with respect to
uniform convergence of gambles: if a sequence of gambles $f_n$ converges
uniformly to a gamble $f$, meaning that for every $\epsilon>0$ there is some
$n_0$ such that $|f_n(x)-f(x)|<\epsilon$ for all $n\geq n_0$ and for all
$x\in\values$, then $\lpr(f_n)$ converges to $\lpr(f)$ and $\upr(f_n)$
converges to $\upr(f)$. In particular, this implies that a coherent lower
prevision defined on some domain $\domain$ can be uniquely extended to a
coherent lower prevision on the uniform closure of $\domain$. As an immediate
corollary, a coherent lower prevision on $\gambles(\values)$ is uniquely
determined by the values it assumes on \emph{simple} gambles, i.e., gambles
that assume only a finite number of values.
\par
We end this section by introducing a number of new notions, which cannot be
found in \citet{walley1991}. They generalise familiar definitions in standard,
measure-theoretic probability to a context where coherent lower previsions are
used as belief models.

\subsection{The distribution of a random variable}
We shall call a subject's coherent lower prevision $\lpr$ on
$\gambles(\values)$, modelling his beliefs about the value that a random
variable $\rv$ assumes in the set $\values$, his \emph{distribution} for that
random variable.
\par
Now consider another set $\avalues$, and a map $\varphi$ from $\values$ to
$\avalues$, then we can consider $\arv:=\varphi(\rv)$ as a random variable
assuming values in $\avalues$. With a gamble $h$ on $\avalues$, there
corresponds a gamble $h\circ\varphi$ on $\values$, whose lower prevision is
$\lpr(h\circ\varphi)$. This leads us to define the distribution of
$\arv=\phi(\rv)$ as the \emph{induced} coherent lower prevision $\alpr$ on
$\gambles(\avalues)$, defined by
\begin{equation*}
  \alpr(h):=\lpr(h\circ\varphi),\quad h\in\gambles(\avalues).
\end{equation*}
For an event $A\subseteq\avalues$, we see that
$I_A\circ\varphi=I_{\varphi^{-1}(A)}$, where
$\varphi^{-1}(A):=\set{x\in\values}{\varphi(x)\in A}$, and consequently
$\alpr(A)=\lpr(\varphi^{-1}(A))$. So we see that the notion of an induced
lower prevision generalises that of an induced probability measure.
\par
Finally, consider a sequence of random variables $\rv_n$, all taking
values in some metric space $S$. Denote by $\conts(S)$ the set of
all continuous gambles on $S$. For each random variable $\rv_n$, we
have a distribution in the form of a coherent lower prevision
$\lpr_{\rv_n}$ on $\gambles(S)$. Then we say that the random
variables \emph{converge in distribution} if for all
$h\in\conts(S)$, the sequence of real numbers $\lpr_{\rv_n}(h)$
converges to some real number, which we denote by $\lpr(h)$. The
limit lower prevision $\lpr$ on $\conts(S)$ that we can define in
this way, is coherent, because a point-wise limit of coherent lower
previsions always is.

\section{Exchangeable random variables}
\label{sec:exchangeable-finite}
We are now ready to recall \citegensec{Section~9.5}{walley1991} notion of
exchangeability in the context of the theory of coherent lower previsions.  We
shall see that it generalises \citeauthor{finetti1937}'s definition for linear
previsions \citep{finetti1937,finetti1975}.

\subsection{Definition and basic properties}\label{sec:basic-definitions}
Consider $N\geq1$ random variables $\rv_1$, \dots, $\rv_N$ taking values in a
non-empty and finite set $\values$.\footnote{We could easily define
  exchangeability for variables that assume values in a set $\values$ that is
  not necessarily finite. But since we only prove interesting results for
  finite $\values$, we have decided to use a finitary context from the
  outset.}  A subject's beliefs about the values that these random variables
$\sample\rv=\tuple{\rv}{N}$ assume jointly in $\values^N$ is given by their
(joint) distribution, which is a coherent lower prevision $\lpr_\values^N$
defined on the set $\gambles(\values^N)$ of all gambles on $\values^N$.
\par
Let us denote by $\permuts_N$ the set of all permutations of $\{1,\dots,N\}$.
With any such permutation~$\pi$ we can associate, by the procedure of lifting,
a permutation of~$\values^N$, also denoted by~$\pi$, that maps any
$\sample{x}=\tuple{x}{N}$ in~$\values^N$ to
$\pi\sample{x}:=(x_{\pi(1)},\dots,x_{\pi(N)})$. Similarly, with any gamble~$f$
on~$\values^N$, we can consider the permuted gamble $\pi f:=f\circ\pi$, or in
other words, $(\pi f)(\sample{x})=f(\pi\sample{x})$ for all
$\sample{x}\in\values^N$.
\par
A subject judges the random variables $\rv_1$, \dots, $\rv_N$ to be
\emph{exchangeable} when he is disposed to exchange any gamble $f$ for the
permuted gamble $\pi f$, meaning that $\lpr_\values^N(\pi
f-f)\geq0$,\footnote{This means that the subject is willing to accept the
  gamble $\pi f-f$, i.e., to exchange $f$ for $\pi f$, in return for any
  positive amount of utility $\epsilon$, however small.} for any permutation
$\pi$. Taking into account the properties of coherence, this means that
\begin{equation*}
  \lpr^N_\values(\pi f-f)=\lpr^N_\values(f-\pi f)=0
\end{equation*}
for all gambles $f$ on $\values^N$ and all permutations $\pi$ in $\permuts_N$.
In this case, we shall also call the joint coherent lower prevision
$\lpr_\values^N$ \emph{exchangeable}. A subject will make an assumption of
exchangeability when there is evidence that the processes generating the
values of the random variables are (physically) similar
\citep[Section~9.5.2]{walley1991}, and consequently the order in which the
variables are observed is not important.
\par
When $\lpr^N_\values$ is in particular a linear prevision $\pr_\values^N$,
exchangeability is equivalent to having $\pr_\values^N(\pi
f)=\pr_\values^N(f)$ for all gambles $f$ and all permutations $\pi$.  Another
equivalent formulation can be given in terms of the (probability) \emph{mass
  function} $p_\values^N$ of $\pr_\values^N$, defined by
\mbox{$p_\values^N(\sample{x}):=\pr_\values^N(\{\sample{x}\})$}.  Indeed, if
we apply linearity to find that $\pr_\values^N(f)
=\sum_{\sample{x}\in\values^N}f(\sample{x})p_\values^N(\sample{x})$, we see
that the exchangeability condition for linear previsions is equivalent to
having \mbox{$p_\values^N(\sample{x})=p_\values^N(\pi\sample{x})$} for all
$\sample{x}$ in $\values^N$, or in other words, the mass function
$p_\values^N$ should be invariant under permutation of the indices.  This is
essentially \citegen{finetti1937} definition for the exchangeability of a
prevision. The following proposition, mentioned by
\citet[Section~9.5]{walley1991}, and whose proof is immediate and therefore
omitted, establishes an even stronger link between Walley's and de Finetti's
notions of exchangeability.

\begin{proposition}\label{prop:exchangeability-lower-envelopes}
  Any coherent lower prevision on $\gambles(\values^N)$ that dominates an
  exchangeable coherent lower prevision, is also exchangeable.  Moreover, let
  $\lpr_\values^N$ be the lower envelope of some set of linear previsions
  $\solp_\values^N$, in the sense that
  \begin{equation*}
    \lpr_\values^N(f)
    =\min\set{\pr_\values^N(f)}{\pr_\values^N\in\solp_\values^N}
  \end{equation*}
  for all gambles $f$ on $\values^N$. Then $\lpr_\values^N$ is exchangeable if
  and only if all the linear previsions~$\pr_\values^N$ in $\solp_\values^N$
  are exchangeable.
\end{proposition}

If a coherent lower prevision $\lpr_\values^N$ is exchangeable, it is
immediately guaranteed to be also \emph{permutable}\footnote{We use the
  terminology in \citet[Section~9.4]{walley1991}.} in the sense that
\begin{equation*}
  \lpr_\values^N(\pi f)=\lpr_\values^N(f)
  \text{ for all gambles $f$ on $\values^N$ and all permutations $\pi$ in
    $\permuts_N$}.
\end{equation*}
The converse does not hold in general. For linear previsions $\pr_\values^N$,
permutability is equivalent to exchangeability, but this equivalence is
generally broken for coherent lower previsions that are not
linear.\footnote{This is an instance of a more general phenomenon: we can
  generally consider two types of invariance of a belief model (a coherent
  lower prevision) with respect to a semigroup of transformations: \emph{weak}
  and \emph{strong} invariance.  The former, of which permutability is a
  special case, tells us that the model or the beliefs are symmetrical
  (symmetry of evidence), whereas the latter, of which exchangeability is a
  special case, reflects that a subject believes there is symmetry (evidence
  of symmetry).  Strong invariance generally implies weak invariance, but the
  two notions in general only coincide for linear previsions.  For more
  details, see \citet{cooman2005c}.}
\par
Clearly, if $\rv_1$, \dots, $\rv_N$ are exchangeable, then any permutation
$\rv_{\pi(1)}$, \dots, $\rv_{\pi(N)}$ is exchangeable as well, and has the
same distribution $\lpr_\values^N$. Moreover, any selection of $1\leq n\leq N$
random variables from amongst the $\rv_1$, \dots, $\rv_N$ are exchangeable
too, and their distribution is given by $\lpr_\values^n$, which is the
$\values^n$-marginal of $\lpr_\values^N$, given by
$\lpr_\values^n(f):=\lpr_\values^N(\cylext{f})$ for all gambles~$f$
on~$\values^n$, where the gamble~$\cylext{f}$ on~$\values^N$ is the
\emph{cylindrical extension} of~$f$ to~$\values^N$, given by
$\cylext{f}\tuple{z}{N}:=f\tuple{z}{n}$ for all $\tuple{z}{N}$ in~$\values^N$.

\begin{running}
  This is the place to introduce our running example. As we go along, we shall
  try to clarify our reasoning by looking at a specific special case, that is
  as simple as possible, namely where the random variables $\rv_k$ we consider
  can assume only two values. So we might be looking at tossing coins, or
  thumbtacks, and consider modelling the exchangeability assessment that the
  order in which these coin flips are considered is of no consequence. More
  generally, our random variables might be the indicators of events:
  $X_k=I_{E_k}$, and then we consider the events $E_1$, \dots, $E_N$ to be
  exchangeable when the order in which they are observed is of no consequence.
  \par
  Formally, we denote the set of possible values for such variables by
  $\booleans=\{0,1\}$, where $1$ and $0$ could stand for heads and tails,
  success and failure, the occurrence or not of an event, and so on. In what
  follows, we shall often call $1$ a success, and $0$ a failure.
  \par
  The joint random variable $\brv=\tuple{\rv}{N}$ then assumes values in the
  space $\booleans^N$, which is made up of all $N$-tuples of zeros and ones.
  As an example, in the case $N=3$, two possible elements of $\booleans^3$ are
  $(1,0,1)$ and $(0,1,1)$. These elements can be related to each other by a
  permutation of the indices, i.e., of the order in which they occur, and
  therefore any exchangeable linear prevision should assign the same
  probability mass to them. And any exchangeable coherent lower prevision is a
  lower envelope of such exchangeable linear previsions. $\lozenge$
\end{running}

\subsection{Count vectors}
Interestingly, exchangeable coherent lower previsions 
have a very simple representation, in terms of sampling without
replacement.\footnote{Actually this is a special case of a much more general
  representation result for coherent lower previsions on a finite space that
  are strongly invariant with respect to a finite group of permutations of
  that space; see \citep{cooman2005c} for more details. Here we give a
  different proof.} To see how this comes about, consider any
$\sample{x}\in\values^N$.  Then the so-called (permutation) \emph{invariant
  atom}
\begin{equation*}
  [\sample{x}]:=\set{\pi\sample{x}}{\pi\in\permuts_N}
\end{equation*}
is the smallest non-empty subset of $\values^N$ that contains $\sample{x}$ and
that is invariant under all permutations $\pi$ in $\permuts_N$. We shall
denote the set of permutation invariant atoms of~$\values^N$
by~$\atoms_\values^N$. It constitutes a partition of the set $\values^N$. We
can characterise these invariant atoms using the \emph{counting maps}
$\cntf_x^N\colon\values^N\to\nats_0$ defined for all $x$ in $\values$ in such
a way that
\begin{equation*}
  \cntf_x^N(\sample{z})=\cntf_x^N\tuple{z}{N}
  :=\abs{\set{k\in\{1,\dots,N\}}{z_k=x}}
\end{equation*}
is the number of components of the $N$-tuple $\sample{z}$ that assume the
value $x$. Here $\abs{A}$ denotes the number of elements in a finite set $A$,
and $\nats_0$ is the set of all non-negative integers (including zero).  We
shall denote by $\bcntf_\values^N$ the vector-valued map from $\values^N$ to
$\nats_0^\values$ whose component maps are the $\cntf_x^N$, $x\in\values$.
Observe that $\bcntf_\values^N$ actually assumes values in the set of
\emph{count vectors}
\begin{equation*}
  \counts_\values^N
  :=\set{\cnt{m}\in\nats_0^\values}{\sum_{x\in\values}m_x=N}.
\end{equation*}
Since permuting the components of a vector leaves the counts invariant,
meaning that $\bcntf_\values^N(\sample{z})=\bcntf_\values^N(\pi\sample{z})$
for all $\sample{z}\in\values^N$ and $\pi\in\permuts_N$, we see that for all
$\sample{y}$ and $\sample{z}$ in $\values^N$
\begin{equation*}
  \sample{y}\in[\sample{z}]
  \iff\bcntf_\values^N(\sample{y})=\bcntf_\values^N(\sample{z}).
\end{equation*}
The counting map $\bcntf_\values^N$ can therefore be interpreted as a
bijection (one-to-one and onto) between the set of invariant atoms
$\atoms_\values^N$ and the set of count vectors $\counts_\values^N$, and we
can identify any invariant atom $[\sample{z}]$ by the count vector
$\cnt{m}=\bcntf_\values^N(\sample{z})$ of any (and therefore all) of its
elements. We shall therefore also denote this atom by $[\cnt{m}]$; and clearly
$\sample{y}\in[\cnt{m}]$ if and only if
$\bcntf_\values^N(\sample{y})=\cnt{m}$. The number of elements $\nu(\cnt{m})$
in any invariant atom $[\cnt{m}]$ is given by the number of different ways in
which the components of any $\sample{z}$ in $[\cnt{m}]$ can be permuted, and
is therefore given by
\begin{equation*}
  \nu(\cnt{m}):=\binom{N}{\cnt{m}}=\dfrac{N!}{\prod_{x\in\values}m_x!}.
\end{equation*}
\par
If the joint random variable $\sample{\rv}=\tuple{\rv}{N}$ assumes the value
$\sample{z}$ in $\values^N$, then the corresponding count vector assumes the
value $\bcntf_\values^N(\sample{z})$ in $\counts_\values^N$. This means that
we can see $\bcntf_\values^N(\sample{\rv})=\bcntf_\values^N\tuple{\rv}{N}$ as
a random variable in $\counts_\values^N$. If the available information about
the values that $\sample{\rv}$ assumes in $\values^N$ is given by the coherent
exchangeable lower prevision $\lpr_\values^N$ --~the distribution of
$\sample{\rv}$~--, then the corresponding uncertainty model for the values
that $\bcntf_\values^N(\sample{\rv})$ assumes in $\counts_\values^N$ is given
by the coherent \emph{induced} lower prevision $\alpr_\values^N$ on
$\gambles(\counts_\values^N)$ --~the distribution of
$\bcntf_\values^N(\sample{\rv})$~--, given by
\begin{equation}\label{eq:counts-distribution}
  \alpr_\values^N(h)
  :=\lpr_\values^N(h\circ\bcntf_\values^N)
  =\lpr_\values^N\bigg(\sum_{\cnt{m}\in\counts_\values^N}
  h(\cnt{m})I_{[\cnt{m}]}\bigg)
\end{equation}
for all gambles $h$ on $\counts_\values^N$. We shall now prove a theorem that
shows that, conversely, any exchangeable coherent lower prevision
$\lpr_\values^N$ is in fact \emph{completely determined} by the corresponding
distribution $\alpr_\values^N$ of the count vectors, also called its
\emph{count distribution}. It also establishes a relationship between
exchangeability and sampling without replacement.
\par
To get where we want, consider an urn with $N$ balls of different types, where
the different types are characterised by the elements $x$ of the set
$\values$. Suppose the \emph{composition} of the urn is given by the count
vector $\cnt{m}\in\counts^N_\values$, meaning that $m_x$ balls are of type
$x$, for $x\in\values$. We are now going to subsequently select (in a random
way) $N$ balls from the urn, without replacing them. Denote by $\arv_k$ the
random variable in $\values$ that is the type of the $k$-th ball selected.
The possible outcomes of this experiment, i.e., the possible values of the
joint random variable $\sample{\arv}=\tuple{\arv}{N}$ are precisely the
elements $\sample{z}$ of the permutation invariant atom $[\cnt{m}]$, and
random selection simply means that each of these outcomes is equally likely.
Since there are $\nu(\cnt{m})$ such possible outcomes, each of them has
probability $1/\nu(\cnt{m})$. Also, any $\sample{z}$ not in $[\cnt{m}]$ has
zero probability of being the outcome of our sampling procedure. This means
that for any gamble $f$ on $\values^N$, its (precise) prevision (or
expectation) is given by
\begin{equation*}
  \muhy_\values^N(f\vert\cnt{m})
  :=\frac{1}{\nu(\cnt{m})}\sum_{\sample{z}\in[\cnt{m}]}f(\sample{z}).
\end{equation*}
The linear prevision $\muhy_\values^N(\cdot\vert\cnt{m})$ is the one
associated with a \emph{multiple hyper-geometric distribution}
\citep[Chapter~39]{johnson1997}, whence the notation. Indeed, for any
$\sample{x}=\tuple{x}{n}$ in $\values^n$, where $1\leq n\leq N$, the
probability of drawing a sequence of balls $\sample{x}$ from an urn with
composition $\cnt{m}$ is given by
\begin{equation*}
  \muhy_\values^N(\{\sample{x}\}\times\values^{N-n}\vert\cnt{m})
  =\frac{\nu(\cnt{m}-\bmu)}{\nu(\cnt{m})}
  =\frac{1}{\nu(\bmu)}\prod_{x\in\values}\binom{m_x}{\mu_x}/\binom{N}{n}
\end{equation*}
where $\bmu=\bcntf_\values^n(\sample{x})$. This means that the probability of
drawing without replacement any sample with count vector $\bmu$ is $\nu(\bmu)$
times this probability [there are that many such samples], and is therefore
given by
\begin{equation*}
  \frac{\nu(\cnt{m}-\bmu)\nu(\bmu)}{\nu(\cnt{m})}
  =\prod_{x\in\values}\binom{m_x}{\mu_x}/\binom{N}{n},
\end{equation*}
which indeed gives the mass function for the multiple hyper-geometric
distribution. For any permutation $\pi$ of $\{1,\dots,N\}$
\begin{equation}\label{eq:muhy-exchangeable}
  \muhy_\values^N(\pi f\vert\cnt{m})
  =\frac{1}{\nu(\cnt{m})}\sum_{\sample{z}\in[\cnt{m}]}f(\pi\sample{z})
  =\frac{1}{\nu(\cnt{m})}\sum_{\pi^{-1}\sample{z}\in[\cnt{m}]}f(\sample{z})
  =\muhy_\values^N(f\vert\cnt{m}),
\end{equation}
since $\pi^{-1}\sample{z}\in[\cnt{m}]$ iff $\sample{z}\in[\cnt{m}]$. This
means that the linear prevision $\muhy_\values^N(\cdot\vert\cnt{m})$ is
exchangeable. The following theorem establishes an even stronger result.

\begin{theorem}[Representation theorem for finite sequences of exchangeable
  variables]\label{theo:finite-representation} Let $N\geq1$ and let
  $\lpr_\values^N$ be a coherent exchangeable lower prevision on
  $\gambles(\values^N)$. Let $f$ be any gamble on $\values^N$.  Then the
  following statements hold:
  \begin{enumerate}[1.]
  \item The gamble $\hat{f}$ on $\values^N$ given by $\hat{f}
    :=\frac{1}{\abs{\permuts_N}}\sum_{\pi\in\permuts_N}\pi f$ is permutation
    invariant, meaning that $\pi\hat{f}=\hat{f}$ for all $\pi\in\permuts_N$.
    It is therefore constant on the permutation invariant atoms of
    $\values^N$, and also given by
    \begin{equation}\label{eq:hatting}
      \hat{f}
      =\sum_{\cnt{m}\in\counts_\values^N}I_{[\cnt{m}]}\muhy_\values^N(f\vert\cnt{m}).
    \end{equation}
  \item $\lpr_\values^N(f-\hat{f})=\lpr_\values^N(\hat{f}-f)=0$, and therefore
    also $\lpr_\values^N(f)=\lpr_\values^N(\hat{f})$.
  \item $\lpr_\values^N(f)=\alpr_\values^N(\muhy_\values^N(f\vert\cdot))$,
    where $\muhy_\values^N(f\vert\cdot)$ is the gamble on $\counts_\values^N$
    that assumes the value $\muhy_\values^N(f\vert\cnt{m})$ in
    $\cnt{m}\in\counts_\values^N$.
  \end{enumerate}
  Consequently a lower prevision on $\gambles(\values^N)$ is exchangeable if
  and only if it has the form $\alpr(\muhy_\values^N(\cdot\vert\cdot))$, where
  $\alpr$ is any coherent lower prevision on $\gambles(\counts_\values^N)$.
\end{theorem}

\begin{proof}
  The first statement is fairly immediate. We therefore turn at once to the
  second statement. Observe that $f-\hat{f}
  =\frac{1}{\abs{\permuts_N}}\sum_{\pi\in\permuts_N}[f-\pi f]$.  Now use the
  coherence [super-additivity and non-negative homogeneity], and the
  exchangeability of the lower prevision $\lpr_\values^N$ to find that
  \begin{equation*}
    \lpr_\values^N(f-\hat{f})
    \geq\frac{1}{\abs{\permuts_N}}
    \sum_{\pi\in\permuts_N}\lpr_\values^N(f-\pi f)
    =0.
  \end{equation*}
  In a completely similar way, we get $\lpr_\values^N(\hat{f}-f)\geq0$.  Since
  it also follows from the coherence [super-additivity] of $\lpr_\values^N$
  that
  $\lpr_\values^N(f-\hat{f})+\lpr_\values^N(\hat{f}-f)\leq\lpr_\values^N(0) =
  0$, we find that indeed
  $\lpr_\values^N(f-\hat{f})=\lpr_\values^N(\hat{f}-f)=0$.  Now let
  $g:=f-\hat{f}$, then $f=\hat{f}+g$ and $\hat{f}=f-g$, and use the coherence
  [super-additivity and accepting sure gains] of $\lpr_\values^N$ to infer
  that
  \begin{equation*}
    \lpr_\values^N(f)
    \geq\lpr_\values^N(\hat{f})+\lpr_\values^N(g)
    =\lpr_\values^N(\hat{f})
    \geq\lpr_\values^N(f)+\lpr_\values^N(-g)
    =\lpr_\values^N(f),
  \end{equation*}
  whence indeed $\lpr_\values^N(f)=\lpr_\values^N(\hat{f})$.
  \par
  To prove the third statement, use
  $\lpr_\values^N(f)=\lpr_\values^N(\hat{f})$ together with
  Equations~\eqref{eq:counts-distribution} and~\eqref{eq:hatting} to find that
  $\lpr_\values^N(f)=\lpr_\values^N(\hat{f})
  =\alpr_\values^N(\muhy_\values^N(f\vert\cdot))$.
  \par
  These statements imply that any exchangeable coherent lower prevision is of
  the form $\alpr(\muhy_\values^N(\cdot\vert\cdot))$, where~$\alpr$ is some
  coherent lower prevision on $\gambles(\counts_\values^N)$. Conversely,
  if~$\alpr$ is any coherent lower prevision on $\gambles(\counts_\values^N)$,
  then $\alpr(\muhy_\values^N(\cdot\vert\cdot))$ is a coherent lower prevision
  on $\gambles(\values^N)$ that is exchangeable: simply observe that for any
  gamble $f$ on $\values^N$ and any $\pi\in\permuts_N$,
  \begin{equation*}
    \alpr(\muhy_\values^N(f-\pi f\vert\cdot))
    =\alpr(\muhy_\values^N(f\vert\cdot)-\muhy_\values^N(\pi f\vert\cdot))
    =\alpr(0)=0,
  \end{equation*}
  taking into account that each $\muhy_\values^N(\cdot\vert\cnt{m})$ is an
  exchangeable linear prevision [Equation~\eqref{eq:muhy-exchangeable}].
\end{proof}

This theorem implies that any exchangeable coherent lower prevision on
$\values^N$ can be associated with, or equivalently, that any collection
of~$N$ exchangeable random variables in~$\values$ can be seen as the result
of,~$N$ random draws without replacement from an urn with~$N$ balls whose
types are characterised by the elements $x$ of $\values$, whose composition
$\cnt{m}$ is unknown, but for which the available information about the
composition is modelled by a coherent lower prevision on
$\gambles(\counts_\values^N)$.\footnote{When $\lpr_\values^N$, and therefore
  also $\alpr_\values^N$, is a linear prevision, i.e., a precise probability
  model, this interpretation follows from the Theorem of Total Probability, by
  interpreting the $\muhy_\values^N(\cdot\vert\cnt{m})$ as conditional
  previsions, and $\alpr_\values^N$ as a marginal. For imprecise models
  $\lpr_\values^N$ and $\alpr_\values^N$, the validity of this interpretation
  follows by analogous reasoning, using Walley's Marginal Extension Theorem;
  see \citet[Section~6.7]{walley1991} and \citet{miranda2006b}. }
\par
That exchangeable linear previsions can be interpreted in terms of sampling
without replacement from an urn with unknown composition, is of course
well-known, and essentially goes back to de Finetti's work on exchangeability;
see \citep{finetti1937} and \citep{cifarelli1996}. \citet{heath1976} give a
simple proof for variables that may assume two values. But we believe our
proof\footnote{\citet[Chapter~9]{walley1991} also mentions this result for
  exchangeable coherent lower previsions.} for the more general case of
exchangeable coherent \emph{lower} previsions and random variables that may
assume more than two values, is conceptually even simpler than Heath and
Sudderth's proof, even though it is a special case of a much more general
representation result \Citep[Theorem~30]{cooman2005c}. The essence of the
present proof in the special case of linear previsions $\pr$ is captured
wonderfully well by \citegensec{Section~3.1}{zabell1992} succinct statement:
``Thus $\pr$ is exchangeable if and only if two sequences having the same
frequency vector have the same probability.''

\begin{running}
  We come back to the simple case considered before, where
  $\values=\booleans$. Any two elements $\sample{x}$ and $\sample{y}$ of
  $\booleans^N$ can be related by some permutation of the indices
  $\{1,\dots,N\}$ iff they have the same number of successes
  $s=\cntf_1^N(\sample{x})=\cntf_1^N(\sample{y})$ (and of course, the same
  number of failures $f=N-s$). We can identify the count space
  $\counts_\booleans^N=\set{(s,f)}{s+f=N}$ with the set
  $\set{s}{s=0,\dots,N}$, and count vectors $\cnt{m}=(s,N-s)$ with the
  corresponding number of successes $s$, which is what we shall do from now
  on.
  \par
  The $2^N$ elements of $\booleans^N$ are divided into $N+1$ invariant atoms
  $[s]$ of elements with the same number of successes $s$, each of which has
  $\nu(s)=\binom{N}{s}=\frac{N!}{s!(N-s)!}$ elements. We have depicted the
  situation for $N=3$ in Figure~\ref{fig:running-1}.
  \par
  \begin{figure}[htb]
    \centering 

\begin{tikzpicture}
  \small
  \tikzstyle{elem}=[color=gray,fill,circle,inner sep=2pt,node distance=1.5cm]
  \tikzstyle{atom}=[pos=.5,below=.8cm]
  \path (-.7,-.5) coordinate (ll)
        (.7,1) coordinate (ur);
  \node[elem,label=above:${(0,0,0)}$] (000) {};
  \node[elem,label=above:${(1,0,0)}$] (100) [right of=000] {};
  \node[elem,label=above:${(0,1,0)}$] (010) [right of=100] {};
  \node[elem,label=above:${(0,0,1)}$] (001) [right of=010] {};
  \node[elem,label=above:${(1,1,0)}$] (110) [right of=001] {};
  \node[elem,label=above:${(1,0,1)}$] (101) [right of=110] {};
  \node[elem,label=above:${(0,1,1)}$] (011) [right of=101] {};
  \node[elem,label=above:${(1,1,1)}$] (111) [right of=011] {};
  \path (000) +(ll) coordinate (0ll)
              +(ur) coordinate (0ur);
  \path (100) +(ll) coordinate (1ll);
  \path (001) +(ur) coordinate (1ur);
  \path (110) +(ll) coordinate (2ll);
  \path (011) +(ur) coordinate (2ur);
  \path (111) +(ll) coordinate (3ll)
              +(ur) coordinate (3ur);
  \draw[rounded corners] (0ll) rectangle (0ur) node[atom] {$s=0$};
  \draw[rounded corners] (1ll) rectangle (1ur) node[atom] {$s=1$};
  \draw[rounded corners] (2ll) rectangle (2ur) node[atom] {$s=2$};
  \draw[rounded corners] (3ll) rectangle (3ur) node[atom] {$s=3$};
\end{tikzpicture}


    \caption{The four invariant atoms $[s]$ in the space
      $\counts_\booleans^3$, characterised by the number of successes $s$.}
    \label{fig:running-1}
  \end{figure}
  \par
  Exchangeability forces each of the elements within an invariant atom $[s]$
  to be `equally likely'. So each $[s]$ is to be considered as a `lump',
  within which probability mass is distributed uniformly. The only freedom
  exchangeability leaves us with, lies in assigning probabilities to the lumps
  $[s]$. This is the essence of Theorem~\ref{theo:finite-representation},
  which tells us that any exchangeable coherent lower prevision
  $\lpr_\booleans^N$ on $\gambles(\booleans^N)$ can be seen as the composition
  of a coherent lower prevision $\alpr_\booleans^N$ on
  $\gambles(\{0,1,\dots,N\})$, representing beliefs about the number of
  successes $s$, and the \emph{hyper-geometric} distributions on $[s]$, which
  guarantee that the probability is distributed uniformly over each of the
  $\nu(s)=\binom{N}{s}$ elements of $[s]$: for any gamble~$f$
  on~$\booleans^N$,
  \begin{equation*}
    \hy^N(f\vert s)
    :=\muhy_\booleans^N(f\vert s,N-s)
    =\frac{1}{\nu(s)}\sum_{\sample{x}\in[s]}f(\sample{x}).\quad\lozenge
  \end{equation*}
\end{running}

For an exchangeable random variable $\brv=\tuple{\rv}{N}$, with (exchangeable)
distribution~$\lpr_\values^N$ on~$\gambles(\values^N)$, we have seen that we
can completely characterise this distribution by the corresponding
distribution of the count vectors $\alpr_\values^N$ on
$\gambles(\counts_\values^N)$.
\par
We have also seen that any selection of $1\leq n\leq N$ random variables from
amongst the \mbox{$\rv_1$, \dots, $\rv_N$} will be exchangeable too, and that
their distribution is given by $\lpr_\values^n$, which is the
$\values^n$-marginal of $\lpr_\values^N$. There is moreover an interesting
relation between the distributions~$\alpr_\values^N$ and $\alpr_\values^n$ of
the corresponding count vectors, which we shall derive in the next section
(Equation~\eqref{eq:time-consistency-counts}).  On the other hand, it is
well-known (see for instance \citet{diaconis1980}; we shall come back to this
in Section~\ref{sec:exchangeable-natex}) that if we have an exchangeable
$N$-tuple $\tuple{\rv}{N}$, it is not always possible to extend it to an
exchangeable $N+1$-tuple. In the next section, we investigate what happens
when we consider exchangeable tuples of arbitrary length.

\section{Exchangeable sequences}
\label{sec:exchangeable-sequences}
\subsection{Definitions}\label{sec:total-joint}
We now generalise the definition of exchangeability from finite to countable
sequences of random variables. Consider a countable sequence $\rv_1$, \dots,
$\rv_n$, \dots\ of random variables taking values in the same non-empty set
$\values$. This sequence is called \emph{exchangeable} if any finite
collection of random variables taken from this sequence is exchangeable. This
is clearly equivalent to requiring that the random variables $\rv_1$, \dots,
$\rv_n$ should be exchangeable for all $n\geq1$.
\par
We can also consider the exchangeable sequence as a single random variable
$\sample{\rv}$ assuming values in the set $\values^\nats$, where $\nats$ is
the set of the natural numbers (positive integers, without zero). Its possible
values $\sample{x}$ are sequences $x_1$, \dots, $x_n$, \dots\ of elements of
$\values$, or in other words, maps from $\nats$ to $\values$. We can model the
available information about the value that $\sample{\rv}$ assumes in
$\values^\nats$ by a coherent lower prevision $\lpr_\values^\nats$ on
$\gambles(\values^\nats)$, called the \emph{distribution} of the exchangeable
random sequence $\sample{\rv}$.
\par
The random sequence $\sample{\rv}$, or its distribution $\lpr_\values^\nats$,
is clearly exchangeable if and only if all its \emph{$\values^n$-marginals}
$\lpr_\values^n$ are exchangeable for $n\geq1$. These marginals
$\lpr_\values^n$ on $\gambles(\values^n)$ are defined as follows: for any
gamble $f$ on $\values^n$,
$\lpr_\values^n(f):=\lpr_\values^\nats(\cylext{f})$, where $\cylext{f}$ is the
cylindrical extension of $f$ to $\values^\nats$, defined by
$\cylext{f}(\sample{x}):=f\tuple{x}{n}$ for all
$\sample{x}=(\ntuple{x}{n},x_{n+1},\dots)$ in~$\values^\nats$. In addition,
the family of exchangeable coherent lower previsions $\lpr_\values^n$,
$n\geq1$, satisfies the following `\emph{time consistency}' requirement:
\begin{equation}\label{eq:time-consistency}
  \lpr_\values^n(f)=\lpr_\values^{n+k}(\cylext{f}),
\end{equation}
for all $n\geq1$, $k\geq0$, and all gambles $f$ on $\values^n$,
where now $\cylext{f}$ denotes the cylindrical extension of $f$ to
$\values^{n+k}$: $\lpr_\values^n$ should be the $\values^n$-marginal
of any $\lpr_\values^{n+k}$.
\par
It follows at once that any finite collection of $n\geq1$ random variables
taken from such an exchangeable sequence has the same distribution as the
first $n$ variables $\rv_1$, \dots, $\rv_n$, which is the exchangeable
coherent lower prevision $\lpr_\values^n$ on $\gambles(\values^n)$.
\par
Conversely, suppose we have a collection of exchangeable coherent lower
previsions $\lpr_\values^n$ on~$\gambles(\values^n)$, $n\geq1$ that satisfy
the time consistency requirement~\eqref{eq:time-consistency}. Then any
coherent lower prevision $\lpr_\values^\nats$ on $\gambles(\values^\nats)$
that has $\values^n$-marginals $\lpr_\values^n$ is exchangeable. The smallest,
or most conservative such (exchangeable) coherent lower prevision is given by
\begin{equation*}
  \lnex_\values^\nats(f)
  :=\sup_{n\in\nats}\lpr_\values^n(\lproj_n(f))
  =\lim_{n\to\infty}\lpr_\values^n(\lproj_n(f)),
\end{equation*}
where $f$ is any gamble on~$\values^\nats$, and its \emph{lower
projection} $\lproj_n(f)$ on~$\values^n$ is the gamble on~$\values^n$
that is defined by
$\lproj_n(f)(\sample{x}):=\inf_{z_k=x_k, k=1,\dots,n}f(\sample{z})$
for all $\sample{x}\in\values^n$, i.e., the lower projection of $f$
on $\sample{x}$ is the infimum of $f$ over the elements of
$\values^\nats$ whose projection on $\values^n$ is $\sample{x}$. See
\citep[Section~5]{cooman2004a} for more details.

\subsection{Time consistency of the count
  distributions}\label{sec:time-consistency}
It will be of crucial interest for what follows to find out what are the
consequences of the time consistency requirement~\eqref{eq:time-consistency}
on the marginals $\lpr_\values^n$ for the corresponding family
$\alpr_\values^n$, $n\geq1$, of distributions of the count vectors
$\bcntf_\values^n\tuple{\rv}{n}$. Consider therefore $n\geq1$, $k\geq0$ and
any gamble $h$ on $\counts_\values^n$. Let $f:=h\circ\bcntf_\values^n$, then
\begin{equation*}
  \alpr_\values^n(h)
  =\lpr_\values^n(f)
  =\lpr_\values^{n+k}(\cylext{f})
  =\alpr_\values^{n+k}(\muhy_\values^{n+k}(\cylext{f}\vert\cdot)),
\end{equation*}
where the first equality follows from Equation~\eqref{eq:counts-distribution},
the second from Equation~\eqref{eq:time-consistency}, and the last from
Theorem~\ref{theo:finite-representation}. Now for any $\cnt{m}'$ in
$\counts_\values^{n+k}$, and any $\sample{z}'=(\sample{z},\sample{y})$ in
$\values^{n+k}=\values^n\times\values^k$ we have that
$\bcntf_\values^{n+k}(\sample{z}')
=\bcntf_\values^n(\sample{z})+\bcntf_\values^k(\sample{y})$ and therefore
\begin{multline}\label{eq:muhy-extend}
  \muhy_\values^{n+k}(\cylext{f}\vert\cnt{m'})\\
  \begin{aligned}
    &=\frac{1}{\nu(\cnt{m}')}
    \sum_{\sample{z}'\in[\cnt{m}']}\cylext{f}(\sample{z}')
    =\frac{1}{\nu(\cnt{m}')}
    \sum_{(\sample{z},\sample{y})\in[\cnt{m}']}f(\sample{z})
    =\frac{1}{\nu(\cnt{m}')}
    \sum_{\substack{\cnt{m}\in\counts_\values^n\\\cnt{m}\leq\cnt{m}'}}
    \sum_{\sample{y}\in[\cnt{m}'-\cnt{m}]}\sum_{\sample{z}\in[\cnt{m}]}
    f(\sample{z})\\
    &=\frac{1}{\nu(\cnt{m}')}
    \sum_{\substack{\cnt{m}\in\counts_\values^n\\\cnt{m}\leq\cnt{m}'}}
    \nu(\cnt{m}'-\cnt{m})\nu(\cnt{m})\muhy_\values^n(f\vert\cnt{m})
    =\sum_{\cnt{m}\in\counts_\values^n}
    \frac{\nu(\cnt{m}'-\cnt{m})\nu(\cnt{m})}{\nu(\cnt{m}')}h(\cnt{m}),
  \end{aligned}
\end{multline}
since $\muhy_\values^n(f\vert\cnt{m})=h(\cnt{m})$, and $\nu(\cnt{m}'-\cnt{m})$
is zero unless $\cnt{m}\leq\cnt{m}'$. So we see that time consistency is
equivalent to
\begin{equation}\label{eq:time-consistency-counts}
  \alpr_\values^n(h)
  =\alpr_\values^{n+k}\bigg(
  \sum_{\cnt{m}\in\counts_\values^n}
  \frac{\nu(\cdot-\cnt{m})\nu(\cnt{m})}{\nu(\cdot)}h(\cnt{m})
  \bigg)
\end{equation}
for all $n\geq1$, $k\geq0$ and $h\in\gambles(\counts_\values^n)$.

\section{A representation theorem for exchangeable sequences}
\label{sec:representation-sequences} \Citet{finetti1937,finetti1975}
has proven a representation result for exchangeable sequences with
linear previsions that generalises
Theorem~\ref{theo:finite-representation}, and where multinomial
distributions take over the r\^ole that the multiple hyper-geometric
ones play for finite collections of exchangeable variables. One
simple and intuitive way \citep[see also][p.~218]{finetti1975} to
understand why the representation result can be thus extended from
finite collections to countable sequences, is based on the fact that
the multinomial distribution can be seen as as limit of multiple
hyper-geometric ones \citep[Chapter~39]{johnson1997}.  This is also
the central idea behind \citegen{heath1976} simple proof of this
representation result in the case of variables that may only assume
two possible values.
\par
However, there is another, arguably even simpler, approach to proving the same
results, which we present here.  It also works for exchangeability in the
context of coherent lower previsions. And as we shall have occasion to explain
further on, it has the additional advantage of clearly indicating what the
`representation' is, and where it is uniquely defined.
\par
We make a start at proving our representation theorem by taking a look at
multinomial processes.

\subsection{Multinomial processes are exchangeable}
Consider a sequence of random variables $\arv_1$, \dots, $\arv_n$, \dots\ that
are mutually independent, and such that each random variable $\arv_n$ has the
same probability mass function $\btheta$: the probability that $\arv_n=x$ is
$\theta_x$ for $x\in\values$.\footnote{In other words, the random variables
  are IID.} Observe that $\btheta$ is an element of the
\emph{$\values$-simplex}
\begin{equation*}
  \simplex_\values=\set{\btheta\in\reals^\values}
  {(\forall x\in\values)(\theta_x\geq0)\text{ and }
    \sum_{x\in\values}\theta_x=1}.
\end{equation*}
Then for any $n\geq1$ and any $\sample{z}$ in $\values^n$ the probability that
$\tuple{\arv}{n}$ is equal to $\sample{z}$ is given by
$\prod_{x\in\values}\theta_x^{\cntf_x(\sample{z})}$, which yields the
\emph{multinomial mass function} \citep[Chapter~35]{johnson1997}. As a result,
we have for any gamble $f$ on $\values^n$ that its corresponding (multinomial)
prevision (expectation) is given by
\begin{align}
  \mult_\values^n(f\vert\btheta) &=\sum_{\sample{z}\in\values^n}f(\sample{z})
  \prod_{x\in\values}\theta_x^{\cntf_x(\sample{z})}
  =\sum_{\cnt{m}\in\counts_\values^n}
  \sum_{\sample{z}\in[\cnt{m}]}f(\sample{z})
  \prod_{x\in\values}\theta_x^{m_x}\notag\\
  &=\sum_{\cnt{m}\in\counts_\values^n}\muhy_\values^n(f\vert\cnt{m})
  \nu(\cnt{m})\prod_{x\in\values}\theta_x^{m_x}\notag\\
  &=\cmult_\values^n(\muhy_\values^n(f\vert\cdot)\vert\btheta),
  \label{eq:multinomial-exchangeable}
\end{align}
where we defined the (count multinomial) linear prevision
$\cmult_\values^n(\cdot\vert\btheta)$ on $\gambles(\counts_\values^n)$ by
\begin{equation}\label{eq:bernstein-cmult-1}
  \cmult_\values^n(g\vert\btheta)
  =\sum_{\cnt{m}\in\counts_\values^n}
  g(\cnt{m})\nu(\cnt{m})\prod_{x\in\values}\theta_x^{m_x},
\end{equation}
where $g$ is any gamble on $\counts_\values^n$.  The corresponding probability
mass for any count vector~$\cnt{m}$, namely\footnote{We assume implicitly that
  $a^0=1$ for all $a\geq0$.}
\begin{equation}\label{eq:bernstein-cmult-2}
  \cmult_\values^n(\{\cnt{m}\}\vert\btheta)
  =\nu(\cnt{m})\prod_{x\in\values}\theta_x^{m_x}
  =:B_{\cnt{m}}(\btheta),
\end{equation}
is the probability of observing some value $\sample{z}$ for $\tuple{\arv}{n}$
whose count vector is $\cnt{m}$. The polynomial function $B_{\cnt{m}}$ on the
$\values$-simplex is called a (multivariate) \emph{Bernstein (basis)
  polynomial}. We have listed a number of very interesting properties for
these special polynomials in the Appendix. One important fact, which we shall
need quite soon, is that the set
$\set{B_{\cnt{m}}}{\cnt{m}\in\counts_\values^n}$ of all Bernstein (basis)
polynomials of fixed degree $n$ forms a basis for the linear space of all
(multivariate) polynomials on $\simplex_\values$ whose degree is at most $n$;
hence their name. If we have a polynomial $p$ of degree $m$, this means that
for any $n\geq m$, $p$ has a \emph{unique} (Bernstein) decomposition
$b^n_p\in\gambles(\counts_\values^n)$ such that
\begin{equation*}
  p=\sum_{\cnt{m}\in\counts_\values^n}b^n_p(\cnt{m})B_{\cnt{m}}.
\end{equation*}
If we combine this with Equations~\eqref{eq:bernstein-cmult-1}
and~\eqref{eq:bernstein-cmult-2}, we find that $b_p^n$ is the unique gamble on~$\counts_\values^n$ such that $\cmult_\values^n(b^n_p\vert\cdot)=p$.
\par
We deduce from Equation~\eqref{eq:multinomial-exchangeable} and
Theorem~\ref{theo:finite-representation} that the linear prevision
$\mult_\values^n(\cdot\vert\btheta)$ on $\gambles(\values^n)$ --~the
distribution of $\tuple{\arv}{n}$~-- is exchangeable, and that
$\cmult_\values^n(\cdot\vert\btheta)$ is the corresponding distribution for
the corresponding count vectors $\bcntf_\values^n\tuple{\arv}{n}$. Therefore
the sequence of IID random variables $\arv_1$, \dots, $\arv_n$, \dots\ is
exchangeable.

\begin{running}
  Let us go back to our example, where $\values=\booleans$. Here the
  $\booleans$-simplex
  $\simplex_\booleans=\set{(\theta,1-\theta)}{\theta\in\unit}$ can be
  identified with the unit interval, and every element
  $\btheta=(\theta,1-\theta)$ can be identified with the probability $\theta$
  of a success.
  \par
  The count multinomial distribution $\cmult_\booleans^n(\cdot\vert\btheta)$
  now of course turns into the (count) \emph{binomial distribution}
  $\cbin^n(\cdot\vert\theta)$ on $\gambles(\{0,\dots,n\})$, given by
  \begin{equation}\label{eq:binomial}
    \cbin^n(g\vert\theta)
    :=\sum_{s=0}^ng(s)\binom{n}{s}\theta^s(1-\theta)^{n-s}
    =\sum_{s=0}^ng(s)B_s^n(\theta)
  \end{equation}
  for any gamble $g$ on the set $\{0,1,\dots,n\}$ of possible values for the
  number of successes~$s$. In this expression, the
  $B_s^n(\theta):=\binom{n}{s}\theta^s(1-\theta)^{n-s}$ are the $n+1$
  (univariate) Bernstein basis polynomials of degree $n$
  \citep{lorentz1986,prautzsch2002}. For fixed $n$, they add up to one and are
  linearly independent, and they form a basis for the linear space of all
  polynomials on $\unit$ of degree at most $n$. $\lozenge$
\end{running}

\subsection{A representation theorem}
Consider the following linear subspace of $\gambles(\simplex_\values)$:
\begin{equation*}
  \poly(\simplex_\values)
  :=\set{\cmult_\values^n(g\vert\cdot)}{n\geq1,g\in\gambles(\counts_\values^n)}
  =\set{\mult_\values^n(f\vert\cdot)}{n\geq1,f\in\gambles(\values^n)},
\end{equation*}
each of whose elements is a \emph{polynomial function} on the
$\values$-simplex:
\begin{align*}
  \cmult_\values^n(g\vert\btheta) &=\sum_{\cnt{m}\in\counts_\values^n}
  g(\cnt{m})\nu(\cnt{m})\prod_{x\in\values}\theta_x^{m_x}
  =\sum_{\cnt{m}\in\counts_\values^n}g(\cnt{m})B_{\cnt{m}}(\btheta),
\end{align*}
and is actually a linear combination of Bernstein basis polynomials
$B_{\cnt{m}}$ with coefficients $g(\cnt{m})$. So $\poly(\simplex_\values)$ is
the linear space spanned by all Bernstein basis polynomials, and is therefore
the set of all polynomials on the $\values$-simplex $\simplex_\values$.
\par
Now if $\yalpr_\values$ is any coherent lower prevision on
$\gambles(\simplex_\values)$, then it is easy to see that the family of
coherent lower previsions $\lpr_\values^n$, $n\geq1$, defined by
\begin{equation}\label{eq:representation}
  \lpr_\values^n(f)=\yalpr_\values(\mult_\values^n(f\vert\cdot)),
  \quad f\in\gambles(\values^n)
\end{equation}
is still exchangeable and time consistent, and the corresponding count
distributions are given by
\begin{equation}\label{eq:representation-counts}
  \alpr_\values^n(f)=\yalpr_\values(\cmult_\values^n(g\vert\cdot)),
  \quad g\in\gambles(\counts_\values^n).
\end{equation}
Here, we are going to show that a converse result also holds: for any time
consistent family of exchangeable coherent lower previsions $\lpr_\values^n$,
$n\geq1$, there is a coherent lower prevision~$\yalpr_\values$ on
$\poly(\simplex_\values)$ such that Equation~\eqref{eq:representation}, or its
reformulation for counts~\eqref{eq:representation-counts}, holds. We shall
call such an $\yalpr_\values$ a \emph{representation}, or representing
coherent lower prevision, for the family $\lpr_\values^n$. Of course, any
representing $\yalpr_\values$, if it exists, is uniquely determined
on~$\poly(\simplex_\values)$.
\par
So consider a family of coherent lower previsions $\alpr_\values^n$ on
$\gambles(\counts_\values^n)$ that are time consistent, meaning that
Equation~\eqref{eq:time-consistency-counts} is satisfied. It suffices to find
an $\yalpr_\values$ such that~\eqref{eq:representation-counts} holds, because
the corresponding exchangeable lower previsions $\lpr_\values^n$ on
$\gambles(\values^n)$ are then uniquely determined by
Theorem~\ref{theo:finite-representation}, and automatically satisfy the
condition~\eqref{eq:representation}.
\par
Our proposal is to \emph{define} the functional $\yalpr_\values$ on the set
$\poly(\simplex_\values)$ as follows: \emph{consider any element $p$ of
  $\poly(\simplex_\values)$. Then, by definition, there is some $n\geq1$ and a
  corresponding unique $b_p^n\in\gambles(\counts_\values^n)$ such that
  $p=\cmult_\values^n(b_p^n\vert\cdot)$.  We then let
  $\yalpr_\values(p):=\alpr_\values^n(b_p^n)$.}
\par
Of course, the first thing to check is whether this definition is consistent:
any polynomial $p$ of degree $m$ has unique representations $b_p^n$ for all
$n\geq m$, which means that we have to check that no inconsistencies can arise
in the sense that
$\alpr_\values^{n_1}(b_p^{n_1})\neq\alpr_\values^{n_2}(b_p^{n_2})$ for some
$n_1,n_2\geq m$. It turns out that this is guaranteed by the \emph{time}
consistency of the $\lpr_\values^n$, or that of the corresponding
$\alpr_\values^n$, as is made apparent by the proof of the following lemma.

\begin{lemma}\label{le:rep-well-defined}
  Consider a polynomial of degree $m$, and let $n_1,n_2\geq m$. Then
  $\alpr_\values^{n_1}(b_p^{n_1})=\alpr_\values^{n_2}(b_p^{n_2})$.
\end{lemma}

\begin{proof}
  We may assume without loss of generality that $n_2\geq n_1$. The Bernstein
  decompositions $b_p^{n_1}$ and $b_p^{n_2}$ are then related by Zhou's
  formula [see Equation~\eqref{eq:bernstein-coefficients-2} in the Appendix]:
  \begin{equation*}
    b_p^{n_2}(\cnt{m}_2)
    =\sum_{\cnt{m}_1\in\counts_\values^{n_1}}
    \frac{\nu(\cnt{m}_2-\cnt{m}_1)\nu(\cnt{m}_1)}{\nu(\cnt{m}_2)}
    b_p^{n_1}(\cnt{m}_1),
    \quad\cnt{m}_2\in\counts_\values^{n_2}.
  \end{equation*}
  Consequently, by the time consistency
  requirement~\eqref{eq:time-consistency-counts}, we indeed get that
  $\alpr_\values^{n_2}(b_p^{n_2})=\alpr_\values^{n_1}(b_p^{n_1})$.
\end{proof}
\noindent
We also have to check whether the functional $\yalpr_\values$ thus defined on
the linear space $\poly_\values$ is a coherent lower prevision. This is
established in the following lemma.

\begin{lemma}
  $\yalpr_\values$ is a coherent lower prevision on the linear space
  $\poly(\simplex_\values)$.
\end{lemma}

\begin{proof}
  We show that $\yalpr_\values$ satisfies the necessary and sufficient
  conditions (P\ref{item:asg})--(P\ref{item:superadd}) for coherence of a
  lower prevision on a linear space.
  \par
  We first prove that (P\ref{item:asg}) is satisfied. Consider any
  $p\in\poly(\simplex_\values)$. Let $m$ be the degree of $p$. We must show
  that $\yalpr_\values(p)\geq\min p$. We find that
  $\yalpr_\values(p)=\alpr_\values^n(b_p^n)\geq\min b_p^n$ for all $n\geq m$,
  because of the coherence [accepting sure gains] of the count lower
  previsions $\alpr_\values^n$.  But Proposition~\ref{prop:ranges} in the
  Appendix tells us that $\min b_p^n\uparrow\min p$, whence indeed
  $\yalpr_\values(p)\geq\min p$.
  \par
  Next, consider any $p$ in $\poly(\simplex_\values)$ and any real
  $\lambda\geq0$.  Consider any $n$ that is not smaller than the degree of
  $p$. Since obviously $b_{\lambda p}^n=\lambda b_p^n$, we get
  \begin{equation*}
    \yalpr_\values(\lambda p)
    =\alpr_\values^n(b_{\lambda p}^n)
    =\alpr_\values^n(\lambda b_{p}^n)
    =\lambda\alpr_\values^n(b_p^n)
    =\lambda\yalpr_\values(p),
  \end{equation*}
  where the third equality follows from the coherence [non-negative
  homogeneity] of the count lower prevision $\alpr_\values^n$. This tells us
  that the lower prevision $\yalpr_\values$ satisfies the non-negative
  homogeneity requirement (P\ref{item:nonneghom}).
  \par
  Finally, consider $p$ and $q$ in $\poly(\simplex_\values)$, and any $n$ that
  is not smaller than the maximum of the degrees of $p$ and $q$. Since
  obviously $b_{p+q}^n=b_p^n+b_q^n$, we get
  \begin{equation*}
    \yalpr_\values(p+q)
    =\alpr_\values^n(b_{p+q}^n)
    =\alpr_\values^n(b_p^n+b_q^n)
    \geq\alpr_\values^n(b_p^n)+\alpr_\values^n(b_q^n)
    =\yalpr_\values(p)+\yalpr_\values(q),
  \end{equation*}
  where the inequality follows from the coherence [super-additivity] of the
  count lower prevision $\alpr_\values^n$. This tells us that the lower
  prevision $\yalpr_\values$ also satisfies the super-additivity requirement
  (P\ref{item:superadd}) and as a consequence it is coherent.
\end{proof}
\noindent
We can summarise the argument above as follows.

\begin{theorem}[Representation theorem for exchangeable
  sequences]\label{theo:representation-sequences} Given a time consistent
  family of exchangeable coherent lower previsions $\lpr_\values^n$ on
  $\gambles(\values^n)$, $n\geq1$, there is a unique coherent lower prevision
  $\yalpr_\values$ on the linear space $\poly(\simplex_\values)$ of all
  polynomial gambles on the $\values$-simplex, such that for all $n\geq1$, all
  $f\in\gambles(\values^n)$ and all $g\in\gambles(\counts_\values^n)$:
  \begin{equation}\label{eq:representation-theorem}
    \lpr_\values^n(f)=\yalpr_\values(\mult_\values^n(f\vert\cdot))
    \quad\text{and}\quad
    \alpr_\values^n(g)=\yalpr_\values(\cmult_\values^n(g\vert\cdot)).
  \end{equation}
\end{theorem}
\noindent
Hence, the belief model governing any countable exchangeable sequence in
$\values$ can be completely characterised by a coherent lower prevision on the
linear space of polynomial gambles on $\simplex_\values$.
\par
In the particular case where we have a time consistent family of exchangeable
\emph{linear} previsions $\pr_\values^n$ on $\gambles(\values^n), n\geq 1$,
then $\yalpr_\values$ will be a linear prevision $\yapr_\values$ on the linear
space $\poly(\simplex_\values)$ of all polynomial gambles on the
$\values$-simplex.  As such, it will be characterised by its values
$\yapr_\values(B_{\cnt{m}})$ on the Bernstein basis polynomials $B_{\cnt{m}}$,
$\cnt{m}\in\counts_\values^n$, $n\geq1$, or on any other basis of
$\poly(\simplex_\values)$.
\par
It is a consequence of coherence that $\yalpr_\values$ is also uniquely
determined on the set $\conts(\simplex_\values)$ of all continuous gambles on
the $\values$-simplex $\simplex_\values$: by the Stone-Weiersta\ss\ theorem,
any such gamble is the uniform limit of some sequence of polynomial gambles,
and coherence implies that the lower prevision of a uniform limit is the limit
of the lower previsions.
\par
This unicity result cannot be extended to more general (discontinuous) types
of gambles: the coherent lower prevision $\yalpr_\values$ is not uniquely
determined on the set of all gambles $\gambles(\simplex_\values)$ on the
simplex: and there may be different coherent lower previsions
$\yalpr_\values^1$ and $\yalpr_\values^2$ on $\gambles(\simplex_\values)$
satisfying Equation~\eqref{eq:representation-theorem}.\footnote{See
  \citet{miranda2006a} for a study of the gambles whose prevision is
  determined by the prevision of the polynomials.} But any such lower
previsions will agree on the class $\poly(\simplex_\values)$ of polynomial
gambles, which is the class of gambles we need in order to characterise the
exchangeable sequence.\footnote{We refrain here from imposing conditions other
  than coherence (e.g., related to $\sigma$-additivity) on such extensions,
  which could guarantee unicity on the set of all measurable gambles; see
  \citet{miranda2006a} for related discussion.}
\par
We now investigate the meaning of the representing lower prevision
$\yalpr_\values$ a bit further. Consider the sequence of so-called
\emph{frequency} random variables
$\bfreqf_n:=\bcntf_\values^n\tuple{\rv}{n}/n$ corresponding to an exchangeable
sequence of random variables $\rv_1$, \dots, $\rv_n$, \dots, and assuming
values in the $\values$-simplex $\simplex_\values$. The distribution
$\lpr_{\bfreqf_n}$ of $\bfreqf_n$, i.e., the coherent lower prevision on
$\gambles(\simplex_\values)$ that models the available information about the
values that $\bfreqf_n$ assumes in $\simplex_\values$, is given by
\begin{equation*}
  \lpr_{\bfreqf_n}(h)
  :=\alpr_\values^n(h\circ\frac{1}{n})
  =\yalpr_\values(\cmult_\values^n(h\circ\frac{1}{n}\vert\cdot)),
  \quad h\in\gambles(\simplex_\values),
\end{equation*}
because we know that $\alpr_\values^n$ is the distribution of\
$\bcntf_\values^n\tuple{\rv}{n}$, and also taking into account
Theorem~\ref{theo:representation-sequences} for the last equality.  Now,
\begin{equation*}
  \cmult_\values^n(h\circ\frac{1}{n}\vert\btheta)
  =\sum_{\cnt{m}\in\counts_\values^n}
  h\big(\frac{\cnt{m}}{n}\big)B_{\cnt{m}}(\btheta)
\end{equation*}
is the \emph{Bernstein approximant} or \emph{approximating Bernstein
  polynomial} of degree $n$ for the gamble $h$, and it is a known result (see
\citep[Section~VII.2]{feller1971b}, \citep[Section~2]{heitzinger2003}) that
the sequence of approximating Bernstein polynomials
$\cmult_\values^n(h\circ\frac{1}{n}\vert\cdot)$ converges uniformly to $h$ for
$n\to\infty$ if $h$ is continuous. So, because $\yalpr_\values$ is defined
uniquely, and is uniformly continuous, on the set $\conts(\simplex_\values)$,
we find the following result, which provides an interpretation for the
representation $\yalpr_\values$, and which can be seen as another
generalisation of de Finetti's Representation Theorem: $\yalpr_\values$ is the
limit of the frequency distributions.

\begin{theorem}\label{theo:convergence-in-distribution}
  For all continuous gambles $h$ on $\simplex_\values$, we have that
  \begin{equation*}
    \lim_{n\to\infty}\lpr_{\bfreqf_n}(h)=\yalpr_\values(h),
  \end{equation*}
  or, in other words, the sequence of distributions $\lpr_{\bfreqf_n}$
  converges point-wise to $\yalpr_\values$ on $\conts(\simplex_\values)$, and
  in this specific sense, \emph{the sample frequencies $\bfreqf_n$ converge in
    distribution}.
\end{theorem}

\begin{running}
  Back to our example, where $\values=\booleans$. Here the Representation
  Theorem (Theorem~\ref{theo:representation-sequences}) states that the
  coherent count lower previsions $\alpr_\booleans^n$, $n\geq1$, for any
  exchangeable sequence of variables in $\booleans$ have the form
  \begin{equation*}
    \alpr_\booleans^n(g)=\yalpr_\booleans(\cbin^n(g\vert\cdot)),
  \end{equation*}
  for all gambles $g$ on the set $\{0,1,\dots,n\}$ of possible numbers of
  successes $s$, where the (count) binomial distribution
  $\cbin^n(\cdot\vert\theta)$ is given by Equation~\eqref{eq:binomial}, and
  $\yalpr_\booleans$ is some coherent lower prevision defined on the set
  $\poly(\unit)$ of all polynomials on $\unit$, which is the set of possible
  values for the probability $\theta$ of a success.
  \par
  This $\yalpr_\booleans$ can be uniquely extended to a coherent lower
  prevision on the set $\conts(\unit)$ of all continuous gambles (functions)
  on $\unit$. And Theorem~\ref{theo:convergence-in-distribution} assures us
  that this $\yalpr_\booleans$ on $\conts(\unit)$ is the `limiting
  distribution' of the frequency of successes $F_1^n=T_1^n\tuple{\rv}{n}/n$,
  as the number of `trials' $n$ goes to infinity.
  \par
  When all the count distributions $\alpr_\booleans^n$ are linear previsions
  $\apr_\booleans^n$, then the representation $\yalpr_\booleans$ is a linear
  prevision $\yapr_\booleans$, and \textit{vice versa}. This linear prevision
  on $\conts(\unit)$, or equivalently, on $\poly(\unit)$ is completely
  determined by (and of course completely determines) its values on any basis
  of the set of polynomials on $\unit$. If we take as a basis the set
  $\set{\theta^n}{n\geq0}$, then we see that $\yapr_\booleans$ is completely
  determined by its (raw) \emph{moment sequence}
  $m_n=\yapr_\booleans(\theta^n)$, $n\geq0$. It is well-known \citep[see for
  instance][Section~VII.3]{feller1971b} that in the case of finitely additive
  probabilities, or linear previsions, a moment sequence uniquely determines a
  distribution function, except in its discontinuity points. And this brings
  us right back to \citegen{finetti1937} version of the Representation
  Theorem: ``la loi de probabilit\'e $\Phi_n(\xi)=P(Y_n\leq\xi)$ tend vers une
  limite pour $n\to\infty$. [\dots] il s'ensuit qu'il existe une loi-limite
  $\Phi(\xi)$ telle que $\lim_{n\to\infty}\Phi_n(\xi)=\Phi(\xi)$ \emph{sauf
    peut-\^etre pour les points de discontinuit\'e}.''\footnote{Our italics.
    In de Finetti's notation, $Y_n$ is our $F_1^n$, and $\Phi_n$ its
    distribution function.} $\lozenge$
\end{running}

\section{Looking at the sample means}
\label{sec:sample-means}
Consider an exchangeable sequence $\rv_1$, \dots, $\rv_n$, \dots, and any
gamble $f$ on $\values$. Then the sequence $f(\rv_1)$, \dots, $f(\rv_n)$,
\dots\ is again an exchangeable sequence of random variables, now taking
values in the finite set $f(\values)$. We are interested in the \emph{sample
  means}
\begin{equation*}
  S_n(f)\tuple{\rv}{n}:=\frac{1}{n}\sum_{k=1}^nf(\rv_k)
\end{equation*}
which form a sequence of random variables in $[\inf f,\sup f]$. For any
$\cnt{m}$ in $\counts_\values^n$ and any $\sample{z}\in[\cnt{m}]$,
\begin{equation*}
  S_n(f)(\sample{z})
  =\frac{1}{n}\sum_{k=1}^nf(z_k)
  =\frac{1}{n}\sum_{x\in\values}m_xf(x)
  =:S_\values\left(f\vert\frac{\cnt{m}}{n}\right)
\end{equation*}
where for each $\btheta\in\simplex_\values$, we have defined the linear
prevision $S_\values(\cdot\vert\btheta)$ on $\gambles(\values)$ by
$S_\values(f\vert\btheta):=\sum_{x\in\values}f(x)\theta_x$.  Observe that
$S_\values(f\vert\cdot)$ is a very special (linear) polynomial gamble on the
$\values$-simplex. We then get
\begin{equation*}
  \muhy_\values^n(S_n(f)\vert\cnt{m})
  =\frac{1}{\nu(\cnt{m})}\sum_{\sample{z}\in[\cnt{m}]}S_n(f)(\sample{z})
  =\frac{1}{\nu(\cnt{m})}\sum_{\sample{z}\in[\cnt{m}]}
  S_\values\left(f\vert\frac{{\cnt{m}}}{n}\right)
  =S_\values\left(f\vert\frac{\cnt{m}}{n}\right)
\end{equation*}
so we find for the distribution $\lpr_{S_n(f)}$ of the sample mean $S_n(f)$,
which is a coherent lower prevision on $\gambles([\inf f,\sup f])$, that
\begin{equation*}
  \lpr_{S_n(f)}(h)
  =\lpr_\values^n(h(S_n(f)))
  =\alpr_\values^n(h(S_\values(f\vert\cdot))\circ\frac{1}{n}),
  \quad h\in\gambles([\inf f,\sup f]).
\end{equation*}
In terms of the representing lower prevision $\yalpr_\values$, we see that
\begin{equation*}
  \cmult_\values^n(h(S_\values(f\vert\cdot)\circ\frac{1}{n})\vert\btheta)
  =\sum_{\cnt{m}\in\counts_\values^n}
  h(S_\values(f\vert\frac{\cnt{m}}{n}))B_{\cnt{m}}(\btheta)
\end{equation*}
is the approximating Bernstein polynomial for the gamble
$h(S_\values(f\vert\cdot))$ on $\simplex_\values$. So for all continuous
gambles $h$ on $[\inf f,\sup f]$, $h(S_\values(f\vert\cdot))$ is a continuous
gamble on $\simplex_\values$, and is therefore the uniform limit of its
sequence of approximating Bernstein polynomials. Since a coherent lower
prevision is uniformly continuous, we see that
\begin{equation}\label{eq:sample-limit-law}
  \lim_{n\to\infty}\lpr_{S_n(f)}(h)=\yalpr_\values(h(S_\values(f\vert\cdot))).
\end{equation}
This tells us that for an exchangeable sequence $\rv_1$, \dots, $\rv_n$,
\dots\ the sequence of sample means $S_n(f)\tuple{\rv}{n}$ converges in
distribution.

\section{Exchangeable natural extension}
\label{sec:exchangeable-natex}
Throughout this paper, we have always considered exchangeable lower previsions
$\lpr_\values^N$ defined on the set $\gambles(\values^N)$ of \emph{all}
gambles on $\values^N$. At first sight, it seems an impossible task to specify
or assess such an exchangeable lower prevision: a subject must specify an
uncountable infinity of supremum acceptable prices, and at the same time keep
track of all the symmetry requirements imposed by exchangeability, as well as
the coherence requirement.
\par
Alternatively, a subject must specify a coherent count lower prevision
$\alpr_\values^N$ on $\gambles(\counts_\values^N)$, and this means specifying
an uncountable infinity of real numbers $\alpr_\values^N(g)$, for all
gambles~$g$ on~$\counts_\values^N$.\footnote{When $\alpr_\values^N$ is a
  linear prevision $\apr_\values^N$, it suffices to specify a finite number of
  real numbers $\apr_\values^N(\{\cnt{m}\})$, for $\cnt{m}$ in
  $\counts_\values^N$, but such an extremely efficient reduction is generally
  not possible for coherent count \emph{lower} previsions~$\alpr_\values^N$.}
\par
Is it therefore realistic, or of any practical relevance, to consider such
exchangeable coherent lower previsions? Indeed it is, and we now want to show
why.

\subsection{The general problem}\label{sec:general}
What will usually happen in practice, is that a subject makes an assessment
that $N$ variables $\rv_1$, \dots, $\rv_N$ taking values in a finite set
$\values$ are exchangeable,\footnote{This is a so-called \emph{structural
    assessment} in \citegen{walley1991} terminology.} and in addition
specifies supremum acceptable buying prices $\lpr(f)$ for all gambles in some
(typically finite, but not necessarily so) set of gambles
$\domain\subseteq\gambles(\values^N)$.  The question then is: \emph{can we
  turn these assessments into an exchangeable coherent lower prevision
  $\lpr_\values^N$ defined on all of $\gambles(\values^\nats)$, that is
  furthermore as small (least-committal, conservative) as possible?}
\par
To answer this question, we begin by looking at the most conservative (i.e.,
point-wise smallest) exchangeable coherent lower prevision
$\lnex_{\permuts_N}$ for $N$ variables.  Since the most conservative coherent
lower prevision on $\gambles(\counts_\values^N)$ is the \emph{vacuous} lower
prevision, given by
$\alpr_\values^N(g)=\min_{\cnt{m}\in\counts_\values^N}g(\cnt{m})$, our
Representation Theorem for finite exchangeable sequences
(Theorem~\ref{theo:finite-representation}) tells us that
\begin{equation}\label{eq:smallest-exchangeable}
  \lnex_{\permuts_N}(f)
  =\min_{\cnt{m}\in\counts_\values^N}\muhy_\values^N(f\vert\cnt{m})
\end{equation}
for all gambles $f$ on $\values^N$, whose corresponding count lower prevision
is vacuous. It models a subject's beliefs about sampling without replacement
from an urn with N balls, where this subject is completely ignorant about the
composition of the urn.
\par
Using this $\lnex_{\permuts_N}$, we can invoke a general theorem we have
proven elsewhere, about the existence of coherent lower previsions that are
(strongly) invariant under a monoid of transformations
\citep[Theorem~16]{cooman2005c} to find
that\footnote{Equation~\eqref{eq:exchangeable-asl} is closely related to the
  avoiding sure loss condition~\eqref{eq:asl}, but where the supremum is
  replaced by the coherent upper prevision $\unex_{\permuts_N}$.  Similarly,
  Equation~\eqref{eq:exchangeable-natex} is related to the
  expression~\eqref{eq:natex} for natural extension, but where the infimum
  operator is replaced by the coherent lower prevision $\lnex_{\permuts_N}$.
  There is a small and easily correctable oversight in the formulation of
  Theorem~16 of \citet{cooman2005c}, as becomes immediately apparent when
  considering its proof: it is there (but should not be) formulated without
  the multipliers $\lambda_k\geq0$.}
\begin{enumerate}[ENE-1.]
\item there are exchangeable coherent lower previsions on
  $\gambles(\values^N)$ that dominate $\lpr$ on $\domain$ if and only if
  \begin{equation}\label{eq:exchangeable-asl}
    \unex_{\permuts_N}\bigg(\sum_{k=1}^n\lambda_k[f_k-\lpr(f_k)]\bigg)\geq0
    \quad\text{ for all $n\geq0$, $\lambda_k\geq0$ and $f_k\in\domain$,
      $k=1,\dots,n$};
  \end{equation}
\item in that case the point-wise smallest (most conservative) exchangeable
  coherent lower prevision $\lnex_{\lpr,\permuts_N}$ on $\gambles(\values^N)$
  that dominates $\lpr$ on $\domain$ is given by
  \begin{equation}\label{eq:exchangeable-natex}
    \lnex_{\lpr,\permuts_N}(f)
    :=\sup\set{\lnex_{\permuts_N}
      \bigg(f-\sum_{k=1}^n\lambda_k[f_k-\lpr(f_k)]\bigg)}
    {n\geq0,\lambda_k\geq0, f_k\in\domain},
  \end{equation}
  and is called the \emph{exchangeable natural extension} of $\lpr$.
\end{enumerate}
If we now combine Equation~\eqref{eq:smallest-exchangeable} with
Equations~\eqref{eq:exchangeable-asl} and~\eqref{eq:exchangeable-natex}, and
define the lower prevision $\alpr$ on the set
\begin{equation*}
  \adomain
  :=\set{\muhy_\values^N(f\vert\cdot)}{f\in\domain}
  \subseteq\gambles(\counts_\values^N)
\end{equation*}
by letting\footnote{Observe that it is necessary that $\alpr(g)$ should be
  finite, in order for the condition~\eqref{eq:exchangeable-asl} to hold.}
\begin{equation*}
  \alpr(g)
  :=\sup\set{\lpr(f)}{\muhy_\values^N(f\vert\cdot)=g,f\in\domain}
\end{equation*}
for all $g\in\adomain$, then it is but a small technical step to prove the
following result.

\begin{theorem}[Exchangeable natural extension]\label{theo:exchangeable-natex}
  There are exchangeable coherent lower previsions on $\gambles(\values^N)$
  that dominate $\lpr$ on $\domain$ if and only if $\alpr$ is a lower
  prevision\footnote{The explicit requirement that $\alpr$ is a lower
    prevision means that $\alpr$ must be nowhere infinite.} on $\adomain$ that
  avoids sure loss. In that case
  $\lnex_{\lpr,\permuts_N}=\lnex_\alpr(\muhy_\values^N(\cdot\vert\cdot))$,
  i.e., the count distribution for the exchangeable natural extension
  $\lnex_{\lpr,\permuts_N}$ of $\lpr$ is the natural extension $\lnex_\alpr$
  of the lower prevision $\alpr$.
\end{theorem}

Since there are quite efficient algorithms \citep{walley2004} for calculating
the natural extension of a lower prevision based on a finite number of
assessments, this theorem not only has intuitive appeal, but it provides us
with an elegant and efficient manner to find the exchangeable natural
extension, i.e., to combine (finitary) local assessments $\lpr$ with the
structural assessment of exchangeability.

\subsection{From \texorpdfstring{$n$}{\textit{n}} to
  \texorpdfstring{$n+k$}{\textit{n}+\textit{k}} exchangeable random
  variables?}
Suppose we have $n$ random variables $\rv_1$, \dots, $\rv_n$, that a subject
judges to be exchangeable, and whose distribution is given by the exchangeable
coherent lower prevision $\lpr_\values^n$ on $\gambles(\values^n)$, with count
distribution~$\alpr_\values^n$ on $\gambles(\counts_\values^n)$. \emph{Can
  this model be extended to a coherent exchangeable model for $n+k$ variables?
  And if so, what is the most conservative such extended model?}
\par
It is well-known that when $\lpr_\values^n$ is a linear prevision, it cannot
generally be extended \citep{diaconis1980}. In the more general case that we
are considering here, we now look at our Theorem~\ref{theo:exchangeable-natex}
to provide us with an elegant answer: the problem considered here is a special
case of the one studied in Section~\ref{sec:general}.
\par
Indeed, if we denote, as before in Section~\ref{sec:total-joint}, by
$\cylext{f}$ the cylindrical extension to $\values^{n+k}$ of the gamble $f$ on
$\values^n$, then we see that the local assessments $\lpr$ are defined on the
set of gambles $\domain:=\set{\cylext{f}}{f\in\gambles(\values^n)}
\subseteq\gambles(\values^{n+k})$ by $\lpr(\cylext{f}):=\lpr_\values^n(f)$,
$f\in\gambles(\values^n)$.  Observe that here $N=n+k$. If we recall
Equation~\eqref{eq:muhy-extend} in Section~\ref{sec:time-consistency}, then we
see that the corresponding set
$\adomain\subseteq\gambles(\counts_\values^{n+k})$ is given by
\begin{equation*}
  \adomain:=\set{\overline{g}}{g\in\gambles(\counts_\values^n)},
\end{equation*}
where for any gamble $g$ on $\counts_\values^n$ and all
$\bmu\in\counts_\values^{n+k}$
\begin{equation*}
  \overline{g}(\bmu)
  :=\sum_{\cnt{m}\in\counts_\values^n}
  \frac{\nu(\cnt{m})\nu(\bmu-\cnt{m})}{\nu(\bmu)}g(\cnt{m})
  =\pr(g\vert\bmu),
\end{equation*}
where $\pr(\cdot\vert\bmu)$ is the linear prevision associated with drawing
$n$ balls without replacement from an urn with composition $\bmu$. Moreover,
for any $h$ in $\adomain$, there is a unique gamble $g$ on $\counts_\values^n$
such that $h=\overline{g}$.\footnote{To see this, consider the polynomial
  $p=\sum_{\bmu\in\counts_\values^{n+k}}h(\bmu)B_\bmu$. Use Zhou's formula
  [Equation~\eqref{eq:bernstein-coefficients-2} in the Appendix] to find that
  if $h=\overline{g}$, then also
  $p=\sum_{\cnt{m}\in\counts_\values^{n}}g(\cnt{m})B_{\cnt{m}}$, and consider
  that expansions in a Bernstein basis are unique.} This implies that the
corresponding lower prevision $\alpr$ on $\adomain$ is given by
\begin{equation*}
  \alpr(\overline{g}):=\alpr_\values^n(g),\quad g\in\gambles(\counts_\values^n).
\end{equation*}
Now observe that
\begin{enumerate}[(a)]
\item $\overline{\lambda}=\lambda$ for all real $\lambda$;
\item $\overline{\lambda g}=\lambda\overline{g}$ for all $g$ in
  $\gambles(\values^n)$ and all real $\lambda$;
\item $\overline{g_1+g_2}=\overline{g}_1+\overline{g_2}$ for all $g_1$ and
  $g_2$ in $\gambles(\values^n)$.
\end{enumerate}
This tells us that $\adomain$ is a linear subspace of
$\gambles(\counts_\values^\nats)$ that contains all constant gambles.
Moreover, because $\alpr_\values^n$ is a coherent lower prevision, we find
that
\begin{enumerate}[(i)]
\item $\alpr(h_1+h_2)\geq\alpr(h_1)+\alpr(h_2)$ for all $h_1$ and $h_2$ in
  $\adomain$;
\item $\alpr(\lambda h)=\lambda\alpr(h)$ for all real $\lambda\geq0$ and all
  $h$ in $\adomain$;
\item $\alpr(h+\lambda)=\alpr(h)+\lambda$ for all real $\lambda$ and all $h$
  in $\adomain$.
\end{enumerate}
Because $\alpr$ and $\adomain$ have these special properties, the condition
for $\lpr_\values^n$ to be extendable to some coherent exchangeable model for
$n+k$ variables, namely that $\alpr$ avoids sure loss on~$\adomain$,
simplifies to $\max\overline{g}\geq\alpr(\overline{g})$ for all
$g\in\gambles(\counts_\values^n)$, i.e., to
\begin{equation*}
  \max_{\bmu\in\counts_\values^{n+k}}
  \sum_{\cnt{m}\in\counts_\values^n}
  \frac{\nu(\cnt{m})\nu(\bmu-\cnt{m})}{\nu(\bmu)}g(\cnt{m})
  \geq\alpr_\values^n(g)
  \quad\text{for all $g\in\gambles(\counts_\values^n)$}.
\end{equation*}
The expression for the natural extension $\lnex_\alpr$ of $\alpr$, applicable
when the above condition holds, can also be simplified significantly, again
because of the special properties of $\alpr$ and~$\adomain$:
\begin{align*}
  \lnex_\alpr(h) &=\sup\set{\inf\Bigl[h
    -\sum_{k=1}^n\lambda_k[\overline{g}_k-\alpr(\overline{g}_k)]\Bigr]}
  {n\geq0,\lambda_k\geq0,g_k\in\gambles(\counts_\values^n)}\\
  &=\sup\set{\inf\left[h-\overline{g}+\alpr(\overline{g})\right]}
  {g\in\gambles(\counts_\values^n)}\\
  &=\sup\set{\alpr(\overline{g}+\inf[h-\overline{g}])}
  {g\in\gambles(\counts_\values^n)}\\
  &=\sup\set{\alpr(\overline{g})}
  {\overline{g}\leq h,g\in\gambles(\counts_\values^n)}\\
  &=\sup\set{\alpr_\values^n(g)}{\overline{g}\leq
    h,g\in\gambles(\counts_\values^n)},
\end{align*}
for all gambles $h$ on $\counts_\values^{n+k}$. The point-wise smallest
extension of $\lpr_\values^n$ to a coherent exchangeable model on
$\gambles(\values^{n+k})$ is then the coherent exchangeable lower prevision
with count distribution $\lnex_\alpr$, because of
Theorem~\ref{theo:exchangeable-natex}.
\par
In the well-known case that $\lpr_\values^n$ is a linear prevision
$\pr_\values^n$, and therefore $\alpr_\values^n$ is also a linear prevision
$\apr_\values^n$, the condition for extendibility can also be written as
\begin{equation*}
  \min_{\bmu\in\counts_\values^{n+k}}\pr(g\vert\bmu)
  \leq\apr_\values^n(g)
  \quad\text{for all $g\in\gambles(\counts_\values^n)$},
\end{equation*}
where on the left hand side we now see the lower prevision of the gamble $g$,
associated with drawing $n$ balls from an urn with $n+k$ balls, of unknown
composition. When this is satisfied, the lower prevision $\alpr$ will actually
be a linear prevision $\apr$ on the linear space $\adomain$, and $\lnex_\apr$
will be the lower envelope of all linear previsions $\apr_\values^{n+k}$ on
$\gambles(\counts_\values^{n+k})$ that extend $\apr$. Similarly, the
exchangeable natural extension will be the lower envelope of all the
exchangeable linear previsions $\pr_\values^{n+k}$ on
$\gambles(\values^{n+k})$ that extend $\pr_\values^n$.

\section{Conclusions}\label{sec:conclusions}
We have shown that the notion of exchangeability has a natural place in the
theory of coherent lower previsions. Indeed, on our approach using Bernstein
polynomials, and gambles rather than events, it seems fairly natural and easy
to derive representation theorems directly for coherent lower previsions, and
to derive the corresponding results for precise probabilities (linear
previsions) as special cases.
\par
Interesting results can also obtained in a context of predictive inference,
where a coherent exchangeable lower prevision for $n+k$ variables is updated
with the information that the first $n$ variables have been observed to assume
certain values. For a fairly detailed discussion of these issues, we refer to
\citet[Section~9.3]{cooman2005c}.
\par
In Section~\ref{sec:sample-means}, we have argued that the sample means
$S_n(f)\tuple{\rv}{n}$ converge in distribution.  It is possible (and quite
easy for that matter) to prove stronger results.  Indeed, using an approach
that is completely similar to the one originally used by \citet{finetti1937},
we can prove that for all non-negative $n$ and $p$:
\begin{equation*}
  \upr_\values^\nats([S_{n+p}(f)-S_{n}(f)]^2)
  \leq 2\frac{p}{n(n+p)}\sup f^2.
\end{equation*}
In other words, for any fixed $p\geq1$, the sequence $S_{n+p}(f)-S_n(f)$
`converges in mean-square' to zero as $n\to\infty$. Even stronger, we find
that for any non-negative $k$ and $\ell$
\begin{equation*}
  \upr_\values^\nats([S_{k}(f)-S_{\ell}(f)]^2)
  \leq 2\frac{\abs{k-\ell}}{k\ell}\sup f^2,
\end{equation*}
and therefore the sequence $S_n(f)$ `Cauchy-converges in mean-square'. These
convergence results can also be used to derive the convergence in distribution
of the $S_n(f)$, but we consider the approach using Bernstein polynomials to
be distinctly more elegant.

\section*{Acknowledgements}
We acknowledge financial support by research grant \mbox{G.0139.01} of the
Flemish Fund for Scientific Research (FWO), and by projects
\mbox{MTM2004-01269}, \mbox{TSI2004-06801-C04-01}. Erik Quaeghebeur's research
was financed by a Ph.D.~grant of the Institute for the Promotion of Innovation
through Science and Technology in Flanders (IWT Vlaanderen).
\par
We would like to thank J\"{u}rgen Garloff for very helpful comments and
pointers to the literature about multivariate Bernstein polynomials.

\appendix
\section{Multivariate Bernstein polynomials}
With any $n\geq0$ and $\cnt{m}\in\counts_\values^n$ there corresponds a
Bernstein (basis) polynomial of degree~$n$ on $\simplex_\values$, given by
$B_{\cnt{m}}(\btheta)=\nu(\cnt{m})\prod_{x\in\values}\theta_x^{m_x}$,
$\btheta\in\simplex_\values$. These polynomials have a number of very
interesting properties \citep[see for instance][Chapters~10
and~11]{prautzsch2002}, which we list here:
\begin{enumerate}[B1.]
\item The set $\set{B_{\cnt{m}}}{\cnt{m}\in\counts_\values^n}$ of all
  Bernstein polynomials of fixed degree $n$ is linearly independent: if
  $\sum_{\cnt{m}\in\counts_\values^n}\lambda_{\cnt{m}}B_{\cnt{m}}=0$, then
  $\lambda_{\cnt{m}}=0$ for all $\cnt{m}$ in
  $\counts_\values^n$.\label{item:bern-linindep}
\item The set $\set{B_{\cnt{m}}}{\cnt{m}\in\counts_\values^n}$ of all
  Bernstein polynomials of fixed degree $n$ forms a partition of unity:
  $\sum_{\cnt{m}\in\counts_\values^n}B_{\cnt{m}}=1$.\label{item:bern-partition}
\item All Bernstein basis polynomials are non-negative, and strictly positive
  in the interior of~$\simplex_\values$.\label{item:bern-positive}
\item The set $\set{B_{\cnt{m}}}{\cnt{m}\in\counts_\values^n}$ of all
  Bernstein polynomials of fixed degree $n$ forms a basis for the linear space
  of all polynomials whose degree is at most $n$.\label{item:bern-basis}
  \setcounter{saveenumi}{\theenumi}
\end{enumerate}
Property B\ref{item:bern-basis} follows from B\ref{item:bern-linindep} and
B\ref{item:bern-partition}. It follows from B\ref{item:bern-basis} that:
\begin{enumerate}[B1.]
  \setcounter{enumi}{\thesaveenumi}
\item Any polynomial $p$ of degree $m$ has a unique expansion in terms of the
  Bernstein basis polynomials of fixed degree $n\geq
  m$,\label{item:bern-expansion}
\end{enumerate}
or in other words, there is a unique gamble $b_p^n$ on $\counts_\values^n$
such that
\begin{equation*}
  p=\sum_{\cnt{m}\in\counts_\values^n}b_p^n(\cnt{m})B_{\cnt{m}}
  =\cmult_\values^n(b_p^n\vert\cdot).
\end{equation*}
This tells us [also use B\ref{item:bern-partition} and
B\ref{item:bern-positive}] that each $p(\btheta)$ is a convex combination of
the Bernstein coefficients $b_p^n(\cnt{m})$, $\cnt{m}\in\counts_\values^n$
whence
\begin{equation}\label{eq:bernstein-coefficients-1}
  \min b_p^n\leq\min p\leq p(\btheta)\leq\max p\leq\max b_p^n.
\end{equation}
It follows from a combination of B\ref{item:bern-partition} and
B\ref{item:bern-basis} that for all $k\geq0$ and all $\bmu$ in
$\counts_\values^{n+k}$,
\begin{equation}\label{eq:bernstein-coefficients-2}
  b_p^{n+k}(\bmu)
  =\sum_{\cnt{m}\in\counts_\values^n}
  \frac{\nu(\cnt{m})\nu(\bmu-\cnt{m})}{\nu(\bmu)}b_p^n(\cnt{m}).
\end{equation}
This is \emph{Zhou's formula} \citep[see][Section~11.9]{prautzsch2002}. Hence
[let $p=1$ and use~B\ref{item:bern-partition}] we find that for all $k\geq0$
and all $\bmu$ in $\counts_\values^{n+k}$,
\begin{equation}\label{eq:bernstein-coefficients-3}
  \sum_{\cnt{m}\in\counts_\values^n}
  \frac{\nu(\cnt{m})\nu(\bmu-\cnt{m})}{\nu(\bmu)}=1.
\end{equation}
The expressions~\eqref{eq:bernstein-coefficients-2}
and~\eqref{eq:bernstein-coefficients-3} also imply that each $b_p^{n+k}(\bmu)$
is a convex combination of the $b_p^n(\cnt{m})$, and therefore $\min
b_p^{n+k}\geq\min b_p^n$ and $\max b_p^{n+k}\leq\max b_p^n$. Combined with
the inequalities in~\eqref{eq:bernstein-coefficients-1}, this leads to:
\begin{equation}\label{eq:bernstein-coefficients-4}
  [\min p,\max p]
  \subseteq[\min b_p^{n+k},\max b_p^{n+k}]
  \subseteq[\min b_p^n,\max b_p^n]
\end{equation}
for all $n\geq m$ and $k\geq0$. This means that the non-decreasing sequence
$\min b_p^n$ converges to some real number not greater than $\min p$, and,
similarly, the non-increasing sequence $\max b_p^n$ converges to some real
number not smaller than $\max p$. The following proposition strengthens this.

\begin{proposition}\label{prop:ranges}
  For any polynomial $p$ on $\simplex_\values$ of degree $m$,
  \begin{equation*}
    \lim_{\substack{n\to\infty\\n\geq m}}[\min b_p^n,\max b_p^n]
    =[\min p,\max p]=p(\simplex_\values).
  \end{equation*}
\end{proposition}

\begin{proof}
  This follows from the fact that the $b_p^{n}$ converge uniformly to the
  polynomial $p$ as $n\to\infty$; see for instance \citet{trump1996}.
  Alternatively, it can be shown
  \citep[see][Section~11.9]{prautzsch2002} that for $n\geq m$
  \begin{equation*}
    b_p^n(\bmu)
    =\sum_{\cnt{m}\in\counts_\values^m}b_p^m(\cnt{m})B_{\cnt{m}}(\frac{\bmu}{n})
    +O(\frac{1}{n})
    =p(\frac{\bmu}{n})+O(\frac{1}{n}),
    \quad\bmu\in\counts_\values^n.
  \end{equation*}
  From this, we deduce that $\min b_p^n \geq \min p+O(\frac{1}{n})$ for any
  $n\geq m$, and as a consequence $\lim_{n\rightarrow \infty, n\geq m}\min
  b_p^n \geq \min p$. If we use now
  Equation~\eqref{eq:bernstein-coefficients-4}, we see that
  $\lim_{n\to\infty,n\geq m}\min b_p^n=\min p$.  The proof of the other
  equality is completely analogous.
\end{proof}


\end{document}